\documentclass[reqno, final, a4paper]{amsart}

\usepackage{amsmath,amssymb,mathtools,xcolor}
\usepackage{latexsym}
\usepackage{subcaption}

\usepackage{amsthm, bm}
\usepackage{url, hyperref}
\usepackage{enumerate}
\usepackage{tikz}

\DeclarePairedDelimiter{\floor}{\lfloor}{\rfloor}
\DeclarePairedDelimiter{\ceil}{\lceil}{\rceil}
\DeclarePairedDelimiter{\card}{|}{|}

\newtheorem{thm}{Theorem}
\newtheorem{prop}[thm]{Proposition}

\newtheorem{lem}[thm]{Lemma}
\newtheorem{cor}[thm]{Corollary}
\numberwithin{thm}{section}
\newtheorem{clm}[thm]{Claim}

\theoremstyle{definition}
\newtheorem{dfn}[thm]{Definition}

\newtheorem{prob}{Problem}

\theoremstyle{remark}

\theoremstyle{definition}

\newcommand{\mb}[1]{\mathbb{#1}}
\newcommand{\mc}[1]{\mathcal{#1}}
\newcommand{\eps}{\varepsilon}
\DeclareMathOperator{\Var}{Var}
\newcommand{\cF}{\mathcal{F}}
\newcommand{\cS}{\mathcal{S}}
\newcommand{\cW}{\mathcal{W}}
\hypersetup{
	colorlinks,
	linkcolor={red!60!black},
	citecolor={green!50!black},
	urlcolor={blue!80!black}
}

\def\COMMENT#1{}
\let\COMMENT=\footnote

\begin{document}

\title[Schur properties of randomly perturbed sets]{Schur properties of randomly perturbed sets}
\author[Shagnik Das]{Shagnik Das}
\address{Department of Mathematics, National Taiwan University}
\email{shagnik@ntu.edu.tw}
\author[Charlotte Knierim] {Charlotte Knierim}
\address{Department of Computer Science, ETH Zurich}
\email{cknierim@inf.ethz.ch}
\author[Patrick Morris]{Patrick Morris}
\address{Department of Mathematics, Universitat Polit\`ecnica de Catalunya}
\email{pmorrismaths@gmail.com}

\thanks{S. Das was supported  by Research supported in part by the Deutsche Forschungsgemeinschaft (DFG, German Research Foundation) Project 415310276 and by the NTU New Faculty Grant 111L7459. P. Morris was supported by the Deutsche Forschungsgemeinschaft (DFG, German Research Foundation) Walter Benjamin program - project number 504502205.}

\begin{abstract}
    
A set $A$ of integers is said to be \emph{Schur} if any two-colouring of $A$ results in monochromatic $x,y$ and $z$ with  $x+y=z$. We study the following problem: how many random integers from $[n]$ need to be added to some $A\subseteq [n]$ to ensure with high probability that the resulting set is Schur? Hu showed in 1980 that when $|A|> \ceil{\tfrac{4n}{5}}$, no random integers are needed, as $A$ is already guaranteed to be Schur. Recently, Aigner-Horev and Person showed that for any dense set of integers $A\subseteq [n]$, adding $\omega(n^{1/3})$ random integers suffices, noting that this is optimal for sets $A$ with $|A|\leq \ceil{\tfrac{n}{2}}$. We close the gap between these two results by showing that if $A\subseteq [n]$ with $|A|=\ceil{\tfrac{n}{2}}+t<\ceil{\tfrac{4n}{5}}$, then adding  $\omega(\min\{n^{1/3},nt^{-1}\})$ random integers will with high probability result in a set that is Schur. Our result is optimal for all $t$, and we further provide a stability result showing that one needs far fewer random integers when $A$ is not close in structure to the extremal examples. We also initiate the study of perturbing sparse sets of integers $A$ by using algorithmic arguments and the theory of hypergraph containers to provide nontrivial upper and lower bounds.

\end{abstract}

\maketitle

\section{Introduction} \label{sec:intro}
A \emph{Schur triple} in a set $A\subseteq \mb{N}$ is a triple  $(x,y,z)\in A^3$ such that $x+y=z$, and we say a set $A \subseteq \mb{N}$ is \emph{$r$-Schur} if any $r$-colouring of the elements in $A$ results in a monochromatic Schur triple.  
Note that the property of $A$ being $1$-Schur is just the property of containing a Schur triple. We call sets that are not $1$-Schur \emph{sum-free}.  This terminology stems from a classic theorem of Schur \cite{schur1916kongruenz} which  asserts that for every $r$, there is some $n_0 = n_0(r)$ such that  $[n]$ is $r$-Schur for all $n \ge n_0$.

Given this, it is natural to ask which subsets of $[n]$ are also $r$-Schur. From an extremal perspective, this leads to the question of establishing the maximum size of a subset $A\subseteq [n]$ that is \emph{not} $r$-Schur. It is a simple exercise to show that if $|A|> \ceil*{\frac{n}{2}}$, $A$ must be $1$-Schur. Taking $A\subseteq[n]$ to be the set of all odd integers or the large integers $\left\{ \floor*{\frac{n}{2}}+1,\ldots,n \right\}$ shows that this is best possible. For $2$-colourings, one can take $A$ to be all integers in $[n]$ that are \emph{not} divisible by $5$, colouring those that are congruent to $1$ or $4$ $(\mbox{mod } 5)$ red and those congruent to $2$ or $3$ $(\mbox{mod } 5)$ blue. This colouring gives no monochromatic Schur triples and hence there exist sets of size $\ceil*{ \frac{4n}{5}}$ that are not $2$-Schur. Hu~\cite{hu1980note} showed with an elegant argument that one can not do better. 
\begin{thm} \label{thm:hu}
For any $n\in \mb{N}$ and $A\subseteq [n]$ with $\card{A}>\ceil*{ \frac{4n}{5}}$, $A$ is $2$-Schur.
\end{thm}
For $r\geq 3$, it remains an open problem to determine what density forces a subset to be $r$-Schur. Abbott and Wang \cite{abbott1977sum} posed this question in 1977 and provided constructions which they conjecture to be best possible, while some upper bounds have been provided in \cite{abbott1977sum,hancock2019independent}.

Deviating from the problem of determining the size of extremal sets, one can also study the behaviour of typical subsets of $[n]$ by adopting a probabilistic perspective. For this, we fix some probability $p=p(n)\in[0,1]$ and randomly sparsify the set $[n]$, defining $[n]_p$ to be the set obtained by taking each integer of $[n]$ into $[n]_p$ independently with probability $p$. The goal is to understand for which $p$ we can expect the resulting set to be $r$-Schur. Here, and throughout, we say an event holds with high probability (whp, for short) if the probability that it holds tends to $1$ as $n$ tends to infinity. Again, establishing the appearance of Schur triples is an easy task and standard tools (the first and second moment methods) give that if $p=o(n^{-2/3})$, then $[n]_p$ is sum-free whp whilst if $p=\omega(n^{-2/3})$ then $[n]_p$ will be $1$-Schur whp. For more colours, the behaviour was determined by Graham, R\"odl and Ruci\'nski \cite{graham1996schur} for $r=2$ and by R\"odl and Ruci\'nski \cite{rodl1997rado} for $r\geq 3$. 

\begin{thm} \label{thm:random}
For any $2\leq r\in \mb{N}$ we have that if $p=o(n^{-1/2})$ then whp $[n]_p$ is not $r$-Schur whilst if $p=\omega(n^{-1/2})$ then whp $[n]_p$ is $r$-Schur.
\end{thm}
For the rest of the paper we restrict to the case $r=2$ and say that a set $A\subseteq [n]$ is \emph{Schur} if it is $2$-Schur.

\subsection*{Randomly perturbed sets of integers}
 The study of randomly perturbed structures appeared with the notion of \emph{smoothed analysis of algorithms}, introduced by Spielman and Teng  \cite{spielman2004smoothed}, where one is interested in interpolating between worst-case and average-case analysis of algorithms by randomly perturbing an input. At a similar time, Bohman, Frieze and Martin~\cite{bfm1} initiated the study of combinatorial properties in randomly perturbed graphs by looking at how many random edges need to be added to an arbitrary dense graph to make it Hamiltonian. As with smoothed analysis, their work bridges the gap between probabilistic and extremal points of view.

This inspired a wealth of results exploring properties of randomly perturbed graphs and hypergraphs. Most pertinent to this work is the study of Ramsey properties. In an analogous fashion to a set of integers being Schur, for $s\in \mb{N}$, we say a graph $G$ is $s$-(edge-)Ramsey if every 2-colouring of the edges of $G$ results in a monochromatic copy of $K_s$. 
A series of results~\cite{das2020ramsey,kst,power} have determined the number of random edges one needs to add to an arbitrary dense graph to ensure that the resulting graph is $s$-Ramsey for all $s\geq 3$. Several variants of this edge-Ramsey problem were also explored in~\cite{das2020ramsey,kst,power} and randomly perturbed graphs have also been studied with respect to vertex-Ramsey~\cite{das2020vertex} and anti-Ramsey~\cite{aigner2019large,aigner2020small} properties. 



Aigner-Horev and Person initiated the study of randomly perturbed structures in the setting of additive combinatorics.
 From our discussion above, if we have a set $A \subseteq [n]$ of integers with $\card{A}\leq \tfrac{4n}{5}$, one can ask how much we need to randomly perturb $A$ in order to obtain a set that is Schur. For dense sets of integers $A$, Aigner-Horev and Person~\cite{aigner2019monochromatic} showed the following. 

\begin{thm} \label{thm:posdense}
Let $\eps > 0$. If $A \subseteq [n]$, $\card{A} \ge \eps n$, and $p = \omega(n^{-2/3})$, then whp $A \cup [n]_p$ is Schur.
\end{thm}

This can be interpreted as saying that any dense set is close to being Schur, since a small random perturbation is enough to force the set to be Schur. From a probabilistic point of view, one can also see that, in comparison to Theorem~\ref{thm:random}, one can save a great deal of randomness by starting with an arbitrary set of positive density.
Note that Theorem~\ref{thm:posdense} is easily seen to be tight for $\card{A} \le \ceil*{\frac{n}{2}}$: taking $A$ to be a sum-free set, we can colour $A$ red and $[n]_p \setminus A$ blue. Then any monochromatic Schur triples must come from $[n]_p$, and the threshold for their appearance, as previously mentioned, is $p = n^{-2/3}$.

\vspace{1mm}

Our first result precisely describes the amount of randomness needed when the size of the deterministic set grows beyond $\tfrac{n}{2}$.

\begin{thm} \label{thm:main_dense}
Let $n$ and $t=t(n)$ be positive integers such that $\ceil*{\frac{n}{2}} + t\leq \ceil*{\tfrac{4n}{5}}$, and define $p(n,t)=  \min \left\{n^{-2/3}, t^{-1} \right\}$. Then the following statements hold.
\begin{enumerate}
\item[(0)]  There exists a set $A \subseteq [n]$ with $\card{A}=\ceil*{\frac{n}{2}} + t$ such that for  $p=o(p(n,t))$, whp $A\cup[n]_p$ is not Schur.
\item[(1)] For all $A \subseteq [n]$ with $\card{A}=\ceil*{\frac{n}{2}} + t$ and $p=\omega(p(n,t))$, whp $A\cup[n]_p$ is Schur.
\end{enumerate}

\end{thm}

In particular, if $|A|\geq \tfrac{n}{2}+\Omega(n)$ then adding a super-constant number of random integers already suffices to force the resulting set to be Schur. Along with Theorems~\ref{thm:hu} and~\ref{thm:posdense}, this completes our understanding of the behaviour of perturbed sets of integers when the starting set is dense, continuing a recent trend in the perturbed setting of exploring the full range of dense starting structures and describing in detail the transition in the random perturbation required (see e.g.~\cite{bottcher2020triangles,bottcher2020cycles,hmt}).

\paragraph{Stability}   Theorem~\ref{thm:main_dense} \emph{(0)} shows that there are sets $A\subseteq [n]$ with $|A|=\ceil*{\frac{n}{2}} + t$ for which (asymptotically) at least $p(n,t)n$ random integers must be added to make the set Schur. Our next result demonstrates that any such set $A$ must have a certain structure. Here we are interested in the case when $t=o(n)$, since when $t=\Omega(n)$, we have $p(n,t)=\Theta(n^{-1})$, and so the question of whether $p(n,t)n$ random integers are necessary simply reduces to determining if $A$ is already Schur or not. When $t=o(n)$, however, we can expect a significant saving for non-extremal examples. In this case, $A$ has size close to $\tfrac{n}{2}$, and a natural candidate for sets requiring many random integers before becoming Schur are those that are close in structure to the extremal sum-free sets. As discussed at the beginning of the introduction, there are two examples of sum-free sets of size $\ceil*{\frac{n}{2}}$, namely the set of odd integers or the set of large integers $\floor*{\frac{n}{2}}+1,\ldots,n$. Moreover, it is well known that there is stability for sum-free sets, in the sense that any large sum-free set must be close in structure to one of these two constructions (see Theorem~\ref{thm:sumstab}). Our next result shows that any set $A$ needing many random integers to be added in order to become Schur must also be close in structure to one of these two examples. 

\begin{thm} \label{thm:stability}
Let $n$ and $t=t(n)$ be positive integers. If $A \subseteq [n]$ with $\card{A} = \ceil*{\tfrac{n}{2}} + t$ and $q=\omega(n^{-1})$ is such that whp $A\cup[n]_q$ is not Schur, then either $\card*{ \left[ \ceil*{\tfrac{n}{2}} , n \right] \setminus A }=O(q^{-1})$ or $A$ contains $O(q^{-2}n^{-1})$ even numbers. 
\end{thm}

Theorem~\ref{thm:stability} shows that even if we only require $\omega(1)$ random integers to be added to $A$ to give a set that is Schur, then we can remove $o(n)$ integers from $A$ to obtain a set contained in one of the two extremal sum-free constructions. Moreover, the dependence of the distance to the sum-free construction on the number of random integers needed is different in the two cases, showing that the set of large numbers is in some sense more sum-free than the set of odd integers. Indeed, note that due to the size constraint, the set $A$ must contain at least $t$ even integers. Thus, if $q=\omega((nt)^{-1/2})$, the second case of Theorem~\ref{thm:stability} cannot occur, and so if $A$ is such that whp $A\cup[n]_q$ is not Schur, then $A$ must be close to the set of large integers. For $t=\omega(n^{1/3})$ we have $(nt)^{-1/2} =o\left( \min \left\{ n^{-2/3}, t^{-1} \right\}\right)$, and so this shows that we can make significant savings in the amount of randomness required by only imposing the condition that $A$ is far from the set of large numbers. 

\paragraph{Sparse base sets}
One can also explore the behaviour of the perturbed model when the base set $A$ is sparse. This direction has recently been explored in the graph  setting~\cite{hahn3513random} and aims to elucidate the full picture of how the randomly perturbed model transitions between the probabilistic and the extremal thresholds. In our setting, the result of Graham, R\"odl and Ruci\'nski (Theorem~\ref{thm:random}) determines the threshold if we have no deterministic elements while  Theorem~\ref{thm:posdense} gives that we can save some randomness when starting from a base set of size $\Omega(n)$. For base sets $A$ of size $o(n)$, we begin by noting that if $|A|=o(n^{1/2})$, we gain nothing compared to starting with an empty base set. Indeed, for any $s=o(n^{1/2})$ we may take $A$ to be $[n]_{q}$, where $q=2sn^{-1}$. Then we have that whp $|A|\geq s$ and, for any $p$, one has $A\cup [n]_p\sim [n]_{p+q-pq}$. By Theorem~\ref{thm:random}, one needs $p=\omega(n^{-1/2})$ in order to ensure the resulting set is Schur whp.

Here, we take a closer look on how many random integers are needed when the size of $A$ transitions from $n^{1/2}$ to $\eps n$. The following theorem provides non-trivial lower and upper bounds on the perturbed threshold in this sparse case. 

\begin{thm}
\label{thm:main_sparse}
Let $n$ and $s=s(n)$ be positive integers with $\Omega\left( n^{1/2} \right) = s\le n/2$. Then the following two statements hold. 
\begin{enumerate}
    \item[(0)] There exists a set $A\subseteq [n]$ with $|A|= s$ such that for $p=o\left((ns)^{-1/3}\right)$, whp $A\cup [n]_p$ is not Schur.
    \item[(1)] For every $A\subseteq [n]$ with $\card{A}=s$ and $p=\omega\left((n^{13}s)^{-1/27}\log n\right)$, whp $A\cup [n]_p$ is Schur.
\end{enumerate}
\end{thm}

\subsection*{Organisation and remarks}

We conclude the introduction with some comments on our proofs and the organisation of the rest of the paper. We will treat the dense base sets of Theorem~\ref{thm:main_dense} and the sparse base sets of Theorem~\ref{thm:main_sparse} separately, as the two settings seem to require very different approaches. Before turning to our proofs, we will outline in Section~\ref{sec:TandT} our notation and collect several number theoretic and probabilistic tools which will be of use to us.

We will then address the dense setting in Section \ref{sec:dense}, where we start by analysing an explicit colouring to prove the $0$-statement of Theorem~\ref{thm:main_dense} in Section~\ref{sec:0dense}. In Section~\ref{sec:1dense} we prove the $1$-statement by adopting the approach of Aigner-Horev and Person~\cite{aigner2019monochromatic}, finding small (11-integer) configurations in our randomly perturbed sets that are themselves Schur. In order to find these configurations, we will use some powerful number theoretic machinery, such as Green's arithmetic removal lemma~\cite{green2005szemeredi}, and our proof will split into cases depending on the structure of our base set $A$. We will also deduce Theorem~\ref{thm:stability} from the proof of the $1$-statement.

In Section~\ref{sec:sparse} we will then turn to the sparse setting, proving Theorem~\ref{thm:main_sparse}. In Section~\ref{sec:0sparse} we prove the $0$-statement, where, in contrast to the dense setting, the proof is non-constructive; we do not give an explicit colouring. Instead, as is common in random Ramsey Theory, we build upon ideas of Graham, R\"odl and Ruci\'nski~\cite{graham1996schur} from their proof of Theorem~\ref{thm:random}, showing that $A\cup [n]_p$ being Schur implies the existence of certain substructures that are whp not present. For the 1-statement, we appeal to the hypergraph container method developed by Saxton and Thomason \cite{saxton2015hypergraph} and independently Balogh, Morris and  Samotij~\cite{BMS15}. In order to use this method in the randomly perturbed setting, we have to apply it to a hypergraph that encodes certain colour configurations that force the appearance of a monochromatic Schur triple when colouring the base set $A$. We believe the use of this hypergraph is one of the most interesting features of our proof and we introduce this, as well as the container method in general, in Section~\ref{sec:containers}. We then use the containers to establish the $1$-statement in Section~\ref{sec:1sparse}.

Finally, in Section~\ref{sec:conc}, we outline some directions for future research. In particular, we discuss the issue of closing the gap between the bounds in Theorem~\ref{thm:main_sparse}, which we find to be a very intriguing open problem.

\section{Terminology and tools} \label{sec:TandT}

Throughout this paper, we will rely on a series of results from number theory and combinatorics. 
For the  convenience of the reader we state the results in this section. First though, we fix some notation and terminology. 

\subsection{Notation} \label{sec:Notation}

As discussed in the introduction, a \emph{Schur triple} in a set $A\subseteq \mathbb{N}$ is a triple $(x,y,z)\in A^3$ such that $x+y=z$. We say that a triple is \emph{degenerate} if $x=y$ and \emph{non-degenerate} otherwise. We say a set $S\subseteq \mathbb{N}$ \emph{hosts} a Schur triple $(x,y,z)$ if $S=\{x\}\cup \{y\} \cup \{z\}$. Note that if a set $S$ hosts a degenerate Schur triple then $|S|=2$ whilst if $S$ hosts a non-degenerate Schur triple then $|S|=3$.  Given a set $A\subseteq \mathbb{N}$, we will sometimes work with the \emph{Schur hypergraph} $\mc H_{\textrm{Schur}}(A)$ generated by $A$, whose vertex set is $A$ and whose edge set consists of all sets contained in $A$ that host Schur triples. 

We say a set  $A\subseteq \mathbb{N}$ is \emph{sum-free} if  it contains no sets that host Schur triples.  We say a set $A\subseteq \mathbb{N}$ is \emph{Schur} if there is no way to partition $A$ into two sum-free sets. In other words, $A$  is Schur if any red/blue-colouring of $A$ results in a monochromatic Schur triple. If $A$ is \emph{not} Schur, we call any red/blue-colouring of $A$ in which both colour  classes form sum-free sets a \emph{Schur colouring}. 

We will work with both tuples (members of $[n]^\ell$ for some $\ell\in \mathbb{N}$) and subsets of $[n]$. We introduce the following notation to ease the exposition. We say a tuple $T\in [n]^\ell$ \emph{contains} a set $S\subseteq [n]$ if all the elements in $S$ appear as entries of $T$. Similarly, we say a set $S\subseteq [n]$ \emph{contains} a tuple $T\in [n]^\ell$ if all the entries of $T$ appear in the set $S$.  We also define the \emph{intersection} of two tuples $T,T'$, denoted $T\cap T'$,  to be the \emph{set} $S\subseteq[n]$ of elements that feature in \emph{both} $T$ and $T'$. Hence $T\cap T'$ is the largest set contained in both $T$ and $T'$. 

Given $p=p(n)$, the random set $[n]_p$ is the set obtained by keeping each element of $[n]=\{1,\ldots,n\}$ independently with probability $p$.

\subsection{Number theoretic tools} \label{sec:NT}

We start with the following arithmetic removal lemma of Green~\cite{green2005szemeredi}.

\begin{thm} \label{thm:removal} 
For every $\eps > 0$ there is a $\delta > 0$ such that if $A \subseteq [n]$ is a set containing at most $\delta n^2$  sets that host Schur triples, then there is a sum-free $A' \subseteq A$ with $\card{A \setminus A'} \le \eps n$.
\end{thm}

The next powerful result we will use is a stability statement for large sum-free sets due to Deshouillers,  Freiman, S\'os and Temkin~\cite{deshouillers1999structure}.

\begin{thm} \label{thm:sumstab}
If $A \subseteq [n]$ is sum-free and $\card{A} > \frac25 n+1$, then either
\begin{itemize}
	\item[(i)] $A$ only consists of odd numbers, or
	\item[(ii)] $\min A > \card{A}$.
\end{itemize}
\end{thm}

We will also often need to find many arithmetic progressions, and the following result of Varnavides \cite{varnavides1959certain} will be repeatedly applied. Here, and throughout, a  \emph{4-AP} in a set $A$ is a sequence $a,a+d,a+2d,a+3d\in A$ and $d\in \mathbb{N}$ is said to be the \emph{difference} of the arithmetic progression.

\begin{thm} \label{thm:APsupersat}
For every $\delta > 0$ there is a $\xi = \xi(\delta) > 0$ such that if $A \subseteq [n]$ is a set with $\card{A} \ge \delta n$, then $A$ contains at least $\xi n^2$ $4$-APs. In particular, there are at least $\xi n$ distinct differences of $4$-APs in $A$.
\end{thm}

\subsection{Probabilistic tools} \label{sec:Prob}

We will use concentration inequalities to guarantee the existence of certain configurations in our random set of integers. 
First, we will use the well-known theorem of Chebyshev, which bounds the deviation from the expectation in terms of the variance; see, for example,~\cite[Chapter 4]{alonspencer}.
\begin{thm}[Chebyshev's inequality]
\label{thm:chebyshev}
Let $X$ be a random variable and let $t> 0$. Then
\[\Pr[\card{X- \mb{E}[X]}\ge t]\le \frac{\Var[X]}{t^2}.\]
\end{thm}

Our second inequality bounds lower tails and can be used to give exponential concentration. 

\begin{thm}[Janson's inequality \cite{janson1990poisson}] 
Let $\Omega$ be a finite set and, for some $m \in \mathbb{N}$, let $S_1,\ldots,S_m\subseteq \Omega$ be a collection of subsets of $\Omega$ (with repetitions allowed). Consider a random subset of $\Omega$ where each element is chosen independently with some probability $p$ and, for $i\in [m]$, let $X_i$ be the indicator random variable for the event that all elements of $S_i$ are chosen in the random set. Let $X=\sum_ {i\in[m]}X_i$ count the number of sets $S_i$ that appear in the random set and, writing $i\sim j$ if $i\neq j$ and $S_i\cap S_j\neq \emptyset$, let
	\[\mu:=\mb{E}[X]=\sum_{i\in [m]}\mb{E}[X_i] \qquad
	\mbox{ and } \qquad \Delta:=\sum_{\substack{i\sim j}}\mb{E}[X_iX_j].\]
	Then for $0\leq t\leq \mu$ we have
	\[\Pr[X\leq \mu -t]\leq e^{-\frac{t^2}{2(\mu+\Delta)}}.\]
	\label{thm:Janson}
\end{thm}

As a key example for how we use Janson's inequality, the following lemma shows that for any large enough collection of Schur triples in $[n]$, our random set whp contains a member of the collection. 

\begin{lem} \label{lem:JansonSTs}
For any $\xi>0$ and $C>0$, there exists a $\zeta>0$ such that the following holds for any $p=p(n)\le Cn^{-1/2}$. If $\cF\subset [n]^3$ is a collection of distinct, non-degenerate Schur triples such that $|\cF|\ge \xi n^2$, then we have that \[\Pr[\cF\cap [n]_p^3= \emptyset]\le e^{-\zeta n^2p^3}.\] 
\end{lem}
\begin{proof}
First, we take a largest subset $\cF'\subseteq \cF$ such that all the sets that host Schur triples in $\cF'$ are \emph{distinct} sets. That is, if $(x,y,z)$ and $(y,x,z)$ both lie in $\cF$, we only take one of these triples in $\cF'$. From this point on, we will only consider triples in $\cF'$ in order to simplify calculations, noting that we have that $\card{\cF'}\geq \frac12 \xi n^2$. For each Schur triple $S\in \cF'$, let $X_S$ be the indicator variable for the event that all three elements of $S$ appear in $[n]_p$. Moreover, let $X=\sum_{S\in \cF'}X_S$. For each $S\in \cF$, we have that $\mb{E}[X_S]=\Pr\left[S\in [n]_p^3\right]=p^3$. Hence, by linearity of expectation, we have that 
\[\mu:=\mb{E}[X]=\sum_{S\in \cF'}\mb{E}[X_S]=\sum_{S\in \cF'}p^3\ge \frac12 \xi n^2p^3.\]
Now for Schur triples $S$ and $S'$, we write $S\sim S'$ if $S\cap S'\ne \emptyset$, recalling our definition of the intersection of tuples from Section~\ref{sec:Notation}.  We will  upper bound \[\Delta:=\sum_{\substack{S\ne S'\in \cF'\\S\sim S'}}\mb{E}[X_SX_{S'}] \quad \mbox{by the larger quantity }\quad \Delta^*:=\sum_{\substack{S\ne S'\in \cS\\S\sim S'}}\mb{E}[X_SX_{S'}],\]
where $\cS$ denotes the set of \emph{all}
non-degenerate Schur triples $(x,y,z) \in [n]^3$ with $x < y$.  We then have  that \[\Delta\le \Delta^*\le 3n^2\cdot 3\cdot 2\cdot p^4+3n^2\cdot 3\cdot 3n\cdot  p^5\le 27(n^2p^4+n^3p^5),\]
where the first summand comes from considering pairs of Schur triples $S\ne S'$ that intersect in 2 elements and the second summand considers pairs that intersect in 1 element. Indeed, in both cases there are at most $3n^2$ choices of $S$ (each pair of elements is contained in at most 3 Schur triples) and then 3 choices of the elements in $S\cap S'$. In the case that $|S\cap S'|=2$, there are then at most 2 choices for $S' \neq S$ given $S\cap S'$. In the case that $|S\cap S'|=1$, there are at most $3n$ choices of Schur triple $S'$ containing the already chosen element of $S\cap S'$.  

Finally we have that $\Delta\le \Delta^*\le 54C^2n^2p^3$ due to our upper bound on $p$, and hence by Theorem~\ref{thm:Janson}, we have that
\[\Pr[X=0]\le \exp \left(-\frac{\mu^2}{2(\mu+\Delta)}\right)\le \exp\left(- \frac{1}{2}\min\left\{\frac{\mu}{2}, \frac{\mu^2}{2\Delta}\right\} \right)\leq \exp(-\zeta n^2p^3), \]
using our lower bound on $\mu$, our upper bound on $\Delta$, and the fact that $\zeta>0$ is chosen small enough with respect to $\xi$ and $C$. 
\end{proof}

We will also use Janson's inequality for  other larger configurations. We make the following definition. 

\begin{dfn} \label{def:wicket}
We say a $9$-tuple \[W=(x_i,y_i,z_i:i=1,2,3)\in [n]^9\] of nine \emph{distinct} elements of $[n]$ is  a \emph{wicket}\footnote{The terminology here is motivated by viewing these configurations in the Schur hypergraph $\mc H_{\textrm{Schur}}([n])$, as well as the sporting interests of one of the authors.} if  $x_i+y_i=z_i$ for $i=1,2,3$ and $x_1+x_2=x_3$. 
\end{dfn}

Note that as we define a wicket to have nine distinct elements, all the four Schur triples in the definition will be non-degenerate. Furthermore, observe that fixing $x_1,x_2,y_1,y_2$ and $y_3$ uniquely determines a wicket. In particular, there are $O(n^5)$ wickets in $[n]$. As with Lemma~\ref{lem:JansonSTs}, we will be interested in showing that if we have a large collection of wickets, a random set will whp contain one of them.

\begin{lem} \label{lem:JansonWickets}
For any $\xi>0$ and $C>1$, there exists a $\zeta>0$ such that the following holds for any $p=p(n)\le Cn^{-1/2}$. If $\cW\subseteq [n]^9$ is a collection of distinct wickets such that $|\cW|\ge \xi n^5$, then we have that \[\Pr[\cW\cap [n]_p^9= \emptyset]\le e^{-\zeta n^5p^9}.\]
\end{lem}

\begin{proof}
As with the proof of Lemma \ref{lem:JansonSTs}, for each wicket $W\in \cW$, we let $X_W$ be the indicator random variable for the event that $W\subseteq [n]_p$. Then, letting $X:=\sum_{W\in \cW}X_W$, we have that \[\mu:=\mathbb{E}[X]=|\cW|p^9\geq \xi n^5p^9.\]
In order to bound $\Delta=\sum\left\{\mathbb{E}[X_WX_{W'}]:W\ne W'\in \cW, W\cap W'\ne \emptyset\right\}$, we consider the number of wickets that can contain a fixed set of a certain size. 

\begin{clm} \label{clm:wicket}
For $\ell=0,1,2,3,4$ and any nonempty set $U\subseteq [n]$ of size at least $8-2\ell$, there are at most $(9!)^{\ell+1} n^{\ell}$ wickets $W\subseteq [n]^9$ containing $U$.
\end{clm}

Before proving the claim, let us see how it implies the lemma. We upper bound $\Delta$ by 
\begin{align*}
    |\cW|\cdot  2^9\cdot p^9 \cdot \Big((9!)^5n^4 p^8 &+(9!)^4n^3  p^7+(9!)^4n^3  p^6+(9!)^3n^2   p^5 \\ &+(9!)^3n^2  p^4+(9!)^2n  p^3+(9!)^2n  p^2 +9!
p +9!\Big).
\end{align*} 
Here we first choose a wicket $W\in \cW$ ($|\cW|$ choices) and then choose some subset of entries of $W$ which will be the intersection $W\cap W'$ (at most $2^9$ choices). The $p^9$ then comes from all the elements of $W$ appearing in $[n]_p$. In the parentheses we then consider the number of choices of $W'$ that intersect $W$ in our already chosen elements of $W\cap W'$. The $i^{th}$ summand in the parentheses corresponds to an intersection of size exactly $i$ with $W$ and hence we have a factor of $p^{9-i}$ to account for the new elements. In each case we use Claim~\ref{clm:wicket} with $\ell = \lceil\frac{8-i}{2}\rceil$ to upper bound the number of wickets $W'$ that can intersect $W$ in our fixed set of size $i$. Note that we require these wickets to intersect in exactly $i$ elements but we can still use the count given by the claim as we are only concerned with an upper bound here. Simplifying, we have that 
\[\Delta \le 2^{104}|\cW|p^9(n^4p^8+n^3p^6+n^2p^4+np^2+1) \le 2^{108}C^8|\cW|p^9, \]
using our upper bound on $p$ in the final inequality. As previously observed, we have $O(n^5)$ wickets in $[n]$. This gives
\[\Delta = O(n^5p^9).\]
Applying Theorem~\ref{thm:Janson}, we have that 
\[\Pr[X=0]\le \exp \left(-\frac{\mu^2}{2(\mu+\Delta)}\right)\le \exp\left(- \frac{1}{2}\min\left\{\frac{\mu}{2}, \frac{\mu^2}{2\Delta}\right\} \right)\le \exp(-\zeta n^5p^9), \]
choosing $\zeta$ sufficiently small. It remains to prove the claim. 

\begin{proof}
We address the cases $\ell=0,1,2,3,4$ in that order. For $\ell=0$, we are interested in a set $U$ of size at least $8$. After a choice of (distinct) label from $\{x_i,y_i,z_i:i=1,2,3\}$ for each element of $U$ (at most $9!$ choices), any wicket $W$ which is labelled to match the labels of $U$ is already fully determined, as the only possible as yet unchosen element of $W$ is contained in a Schur triple with two already labelled elements in $U$. 

For $\ell=1$, consider a set $U$ of size at least 6 and choose distinct labels for the elements of $U$ from $\{x_i,y_i,z_i:i=1,2,3\}$ (at most $9!$ choices of labels). Now consider the number of wickets $W$ which contain $U$ and whose labels coincide with how we have labelled $U$. Note that there are two indices $i_1,i_2\in [3]$, such that for $i=i_1,i_2$, at least two of the labels in the set  $\{x_i,y_i,z_i\}$ have been assigned to $U$ already. As we must have that $x_i+y_i=z_i$ for $i=i_1,i_2$ and any wicket whose labels are compatible with the labels of $U$, we have that the labels  $\{x_i,y_i,z_i:i=i_1,i_2\}$ of any such wicket are already determined. Hence, letting $j=[3]\setminus\{i_1,i_2\}$, we have that $x_j$ is also determined. A choice of $y_j$ (at most $n$ choices) then completely determines $W$. 

For $\ell=2$, consider a set $U$ of size at least $4$ and choose labels for the elements of $U$ from  $\{x_i,y_i,z_i:i=1,2,3\}$ (at most $9!$ choices). Note that there is some $i^*\in [3]$ such that at least two of the labels $\{x_{i^*},y_{i^*},z_{i^*}\}$ have been placed on elements of $U$. Let $j^*\in [3]\setminus \{i^*\}$ be such that at least one of the labels of $\{x_{j^*},y_{j^*},z_{j^*}\}$ have been placed on $U$ (note that such a $j^*$ is guaranteed to exist). If at least two of the labels of $\{x_{j^*},y_{j^*},z_{j^*}\}$ appear on $U$, then any wicket whose labelling is compatible with $U$ must contain the 6 elements $\{x_k,y_k,z_k:k=i^*,j^*\}$ and so we can appeal to the upper bound for the $\ell=1$ case. Similarly, if only one of the labels $\{x_{j^*},y_{j^*},z_{j^*}\}$ appears on $U$, then a choice of a further element in $[n]$ (at most $n$ choices)  labelled with another label from $\{x_{j^*},y_{j^*},z_{j^*}\}$ (the first free label according to the predetermined ordering of the tuple) gives a set of size 6 that any wicket compatible with the already labelled elements of $[n]$ must contain. Again, we can appeal to induction in this case and conclude that there are at most $9!\cdot (9!)^2n\cdot n\le ( 9!)^3n^2$ wickets containing $U$, as required. 

The case $\ell=3$ is similar. We consider a set $U\subseteq [n]$ of size at least $2$. If $|U|\ge 4$ then we are done by the $\ell=2$ case. Hence we can assume that $|U|=2$ or $|U|=3$. We choose labels for the elements of $U$ from  $\{x_i,y_i,z_i:i=1,2,3\}$ (fewer than $9!$ choices) and consider first the case that there is some $i_0\in [3]$ such that exactly one label of $\{x_{i_0},y_{i_0},z_{i_0}\}$ has been assigned. By making a choice of an element to receive one of the other labels in $\{x_{i_0},y_{i_0},z_{i_0}\}$ (at most $n$ choices) we fix all the elements in $\{x_{i_0},y_{i_0},z_{i_0}\}$, and we therefore have at least $4$ elements already labelled in $[n]$. Counting wickets containing these $4$ elements reduces to the $\ell=2$ case and we can use induction. Similarly, if there is no such $i_0$, then there must be some $i_0'$ such that at least  two of the labels in $\{x_{i'_0},y_{i'_0},z_{i'_0}\}$ have been used in $U$ and so all the elements in $\{x_{i'_0},y_{i'_0},z_{i'_0}\}$  are determined for any wicket containing $U$. Choosing any further element (at most $n$ choices) and labelling it (with the first free label of the tuple) gives 4 fixed elements and we reduce again to the $\ell=2$ case. 

Finally when $\ell=4$, we consider sets $U$ with $|U|\geq 1$ (as $U$ is nonempty by assumption). If $|U|\ge 2$, we are already done so we can assume $ |U|=1$. Choose a label for the element of $U$ (at most 9 choices) and choose another element of $[n]$ to label with the first free label. Counting wickets containing these two elements reduces to the $\ell=3$ case and we have at most $9\cdot n \cdot (9!)^4 n^3\le (9!)^5 n^4$ wickets containing $U$ as required. 
\end{proof} \end{proof}

\section{Dense base sets} 
\label{sec:dense}
In this section we will prove Theorems~\ref{thm:main_dense} and~\ref{thm:stability}.
Let us first look at the $0$-statement of Theorem~\ref{thm:main_dense}: we will prove the lower bound by analysing a particular colouring of an explicit dense set that has been randomly perturbed. 

\subsection{Proof of the $0$-statement of Theorem~\ref{thm:main_dense}}
\label{sec:0dense}

When $t=\Omega(n)$, and $\card{A} = \ceil*{\frac{n}{2}} + t\leq \ceil*{\tfrac{4n}{5}}$, we simply take $A$ to be any set which is not Schur (for example, the construction that removes the integers divisible by $5$, as discussed in the introduction). Then for $p=o(n^{-1})$, an application of Markov's inequality gives that whp $[n]_p$ is empty and $A\cup [n]_p$ remains 2-colourable without monochromatic Schur triples. 

\begin{figure}[h!]
    \centering
    \includegraphics[scale= 0.78]{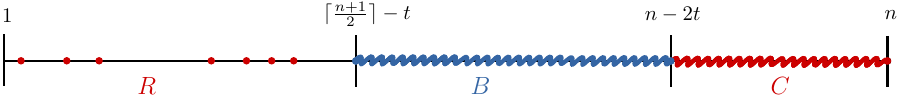}
    \caption{Visualisation of the lower bound construction}
    \label{fig:lower_bound}
\end{figure}

For $1\leq t=o(n)$, let $A = [\ceil*{\frac{n+1}{2}} - t, n]$ and let $p = o \left( \min \{ n^{-2/3}, t^{-1} \} \right)$. Write $B = [\ceil*{\frac{n+1}{2}} - t, n-2t]$, $C = [n-2t+1, n]$, and $R = [n]_p \setminus A = [\floor{\frac{n}{2}} - t ]_p$, noting that $A \cup [n]_p = B \cup C \cup R$. We colour $B$ blue and $C \cup R$ red, as pictured in Figure~\ref{fig:lower_bound}. Note that $B$ is sum-free, and therefore we have no monochromatic Schur triples in blue. We also have that $C$ is sum-free, and since $\min C > 2 \max R$, the only possible monochromatic red Schur triples are of the form $x + y = z$ with $x,y,z \in R$ or with $x \in R$ and $y, z \in C$. The former amounts to the random set containing a Schur triple, which we know whp does not happen for $p = o(n^{-2/3})$. For the latter, we require the element $x$ to belong to the difference set $C - C$. Since $C$ is an interval of length $2t$, there are $2t - 1$ possible differences. As $p = o(t^{-1})$, whp none of these elements $x$ appear in $R$. Thus, this colouring has no monochromatic Schur triples whp, thereby demonstrating that $A \cup [n]_p$ is whp not Schur.\qed

\subsection{Proof of the $1$-statement of Theorem~\ref{thm:main_dense}}

\label{sec:1dense}

We use the following variations of a fact used by Aigner-Horev and Person~\cite{aigner2019monochromatic}, which observes that certain sets are Schur. Our proof of the $1$-statement of Theorem~\ref{thm:main_dense} will then reduce to proving the existence of one of these sets in the randomly perturbed set.

\begin{prop} \label{prop:11set}
Let $a,x,d\in [n]$. Then the following two sets are Schur:
\begin{itemize}
	\item[(i)] $L_1(a,x,d) = \{ d, x, x+d, a, a+d, a+2d, a+3d, a+x, a+x+d, a+x+2d, a+x+3d \}$, and
	\item[(ii)] $L_2(a,x,d) = \{ d, x-d, x, a, a+d, a+2d, a+3d, x-a-3d, x-a-2d, x-a-d, x-a \}$.
\end{itemize}
\end{prop}

\begin{figure}
    \centering
\begin{minipage}{0.43\linewidth}
    \centering
\includegraphics[scale=0.33]{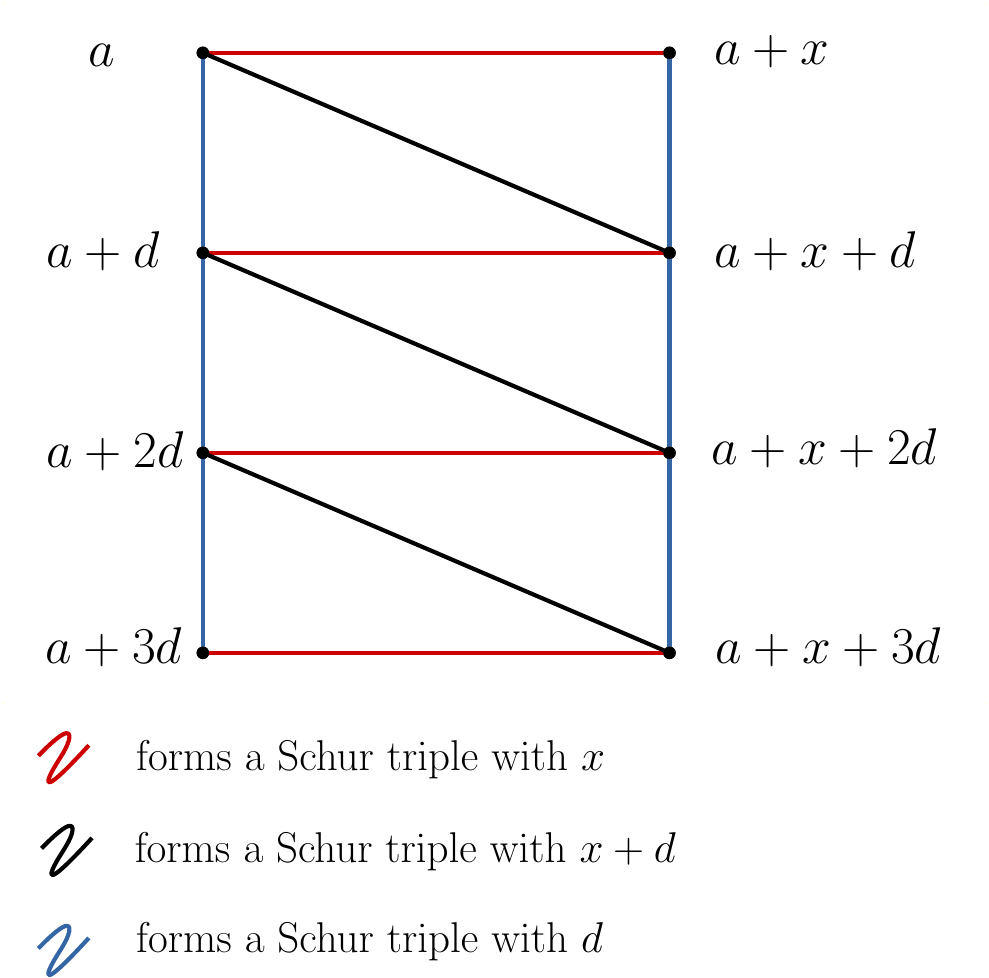}
\label{fig:L1}
\subcaption{A visualisation of $L_1(a,x,d)$}
\end{minipage}
\vspace{0.01\linewidth}
\begin{minipage}{0.43\linewidth}
    \centering
\includegraphics[scale=0.33]{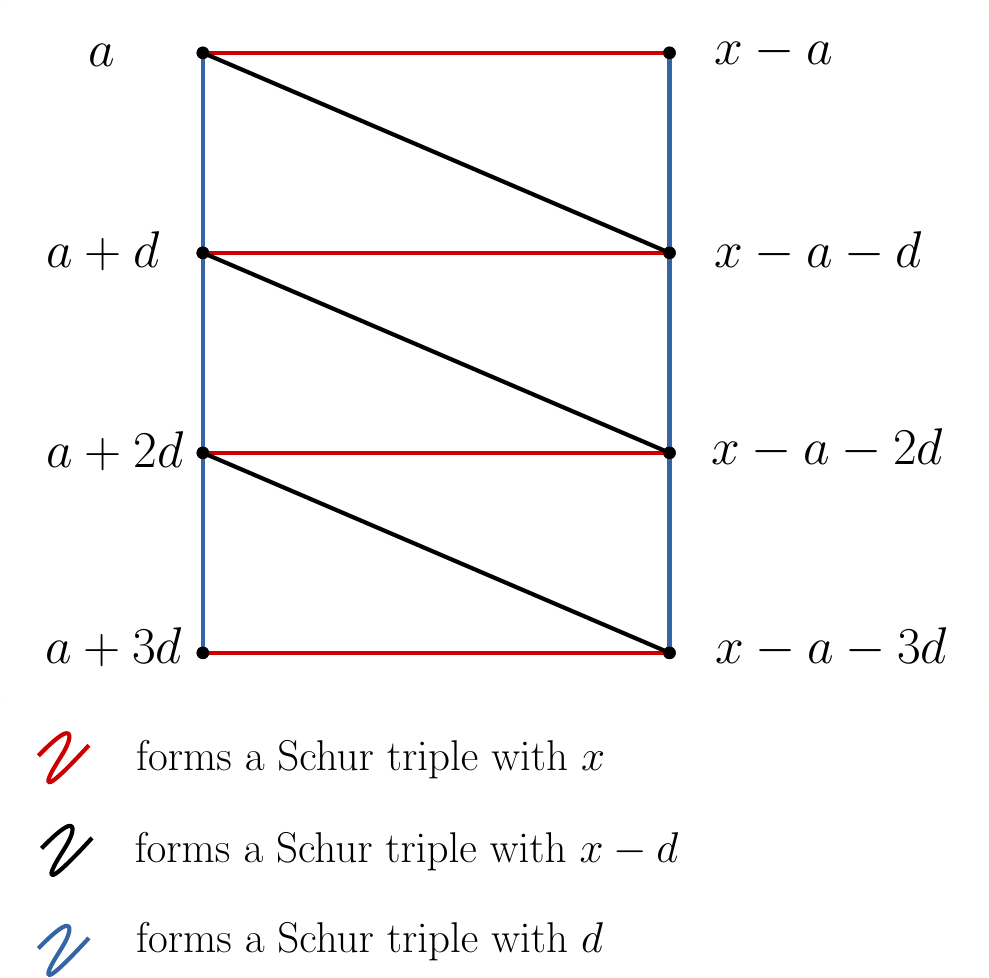}
\label{fig:L2}
\subcaption{A visualisation of $L_2(a,x,d)$}
\end{minipage}
\caption{A comparison of the sum structure of $L_1(a,x,d)$ and $L_2(a,x,d)$.}
\label{fig:compare_struct}
\end{figure}

The proof of this proposition 
follows from a simple case analysis and we omit the details. One can also derive it from the proof of Lemma 2 in \cite{aigner2019monochromatic}. Indeed, Aigner-Horev and Person define a configuration similar\footnote{For any $a,x,d$ we have that $L_1(a,x,d) \supseteq \mc L(a+d,x+a+2d,d)$ where $\mc L$ is as defined in \cite{aigner2019monochromatic}.} to our $L_1(a,x,d)$ and prove that such a configuration is Schur. The proof can be followed directly to prove that $L_1(a,x,d)$ is Schur for all $a,x,d\in [n]$. Moreover, the proof relies solely on the Schur triples depicted in Figure~\ref{fig:compare_struct} and, as shown in the figure, there is an isomorphism between these Schur triples in $L_1(a,x,d)$ and the Schur triples in $L_2(a,x,d)$, thus verifying that $L_2(a,x,d)$ is also Schur.

Given an element $x \in A$, we define $S^+_A(x) = \{ y \in A: x + y \in A \}$, $S^-_A(x) = \{ y \in A : x - y \in A \}$, and $S_A(x) = S^+_A(x) \cup S^-_A(x)$. The following result shows that it will suffice to find some structure in these sets of candidates for Schur triples in the set $A$.

\begin{lem} \label{lem:chebyshev}
Suppose $\lambda=\lambda(n),\kappa=\kappa(n)$ are integers with $\kappa\geq 2$ and we have a set $X \subseteq A \subseteq  [n]$ of size $\lambda$, and that, for each $x \in X$, there is a set $D_x$ of size $\kappa$ such that for every $d \in D_x$, either $S^+_A(x)$ or $S^-_A(x)$ contains a $4$-AP with common difference $d$. If $p = \omega \left( \max \{ ( \lambda \kappa)^{-1/2}, \kappa^{-1} \} \right)$, then $A \cup [n]_p$ is Schur whp.
\end{lem}

Before proving this lemma, let us see how it implies the $1$-statement of Theorem~\ref{thm:main_dense}. 

\begin{proof}[Proof of the $1$-statement of Theorem~\ref{thm:main_dense}]
First observe that if $t = O(n^{2/3})$, then $p(n,t) = n^{-2/3}$, and it follows from Theorem~\ref{thm:posdense} that $A \cup [n]_p$ is whp Schur. Hence, we may assume that $t = \omega (n^{2/3})$, and that $p(n,t) = t^{-1}$.

We split the proof into two cases, depending on the number of Schur triples in $A$. Set $\eps = \tfrac{1}{50}$ and let $\delta = \delta_{\ref{thm:removal}}(\tfrac{1}{50}) > 0$ be the resulting value from Theorem~\ref{thm:removal}. Recall that $|A|= \ceil{\frac{n}{2}}+t$, and note that by monotonicity we may assume $t \le \eps n$. Indeed, if this is not the case then we can shrink $A$ to a subset of size $\ceil{\frac{n}{2}}+\eps n$ and work with this base set instead. 

\paragraph{Case I: there are at least $\delta n^2$ Schur triples in $A$}

Let $X = \{ x \in A : \card{S^+_A(x)} \ge \tfrac12 \delta n \}$. By counting Schur triples, we have
\[ \delta n^2 \le \sum_{x \in A} \card{S^+_A(x)} \le n \cdot \tfrac12 \delta n + \card{X} \cdot n, \]
and so $\card{X} \ge \tfrac12 \delta n$.

Now, by Theorem~\ref{thm:APsupersat}, there is some $\xi > 0$ such that, for each $x \in X$, there is a set $D_x$ of at least $\xi n$ values $d$ such that $S^+_A(x)$ contains a $4$-AP with common difference $d$.
We may therefore apply Lemma~\ref{lem:chebyshev} with $\lambda = \ceil*{\tfrac12 \delta n}$ and $\kappa = \ceil*{\xi n}$ to deduce that, when $p = \omega (n^{-1})$, $A \cup [n]_p$ is Schur whp. 

\paragraph{Case II: there are fewer than $\delta n^2$ Schur triples in $A$}

By Theorem~\ref{thm:removal}, we can remove at most $\eps n$ elements from $A$ to obtain a sum-free subset $A' \subseteq [n]$. It follows that $\card{A'} \ge (\tfrac12 - \eps)n$, and hence we can apply Theorem~\ref{thm:sumstab} to obtain structural information about $A'$ --- it either consists entirely of odd integers or of large integers.

\paragraph{Case II.1: $A'$ is contained in the odd integers}

Since $\card{A'} \ge ( \tfrac12 - \eps )n$, it follows that at most $\eps n$ odd integers are missing from $A$. Furthermore, letting $X\subseteq A$ be the set of even integers in $A$, we have that $\ell:=\card{X}\geq t$, since $\card{A} = \ceil{\frac{n}{2}} + t$.
Now for $x\in X$, if $x \le \frac{n}{2}$, there are at least $\frac{n}{4}$ pairs $a,b\in [n]$ with $a, b$ odd such that $a + x = b$. Since each missing odd integer appears in at most two of these pairs, it follows that $\card{S^+_A(x)} \ge (\tfrac14 - 2 \eps)n \ge \tfrac18 n$. Thus, by Theorem~\ref{thm:APsupersat}, there is some $\xi > 0$ such that there is a set $D_x$ of at least $\xi n$ distinct differences of $4$-APs contained in $S^+_A(x)$.

On the other hand, if $x > \frac{n}{2}$, then there are at least $\frac{n}{4}$ pairs $a,b\in [n]$ with $a, b$ odd such that $a+b=x$. Each missing odd integer appears in at most one such pair, and so $\card{S^-_A(x)} \ge (\tfrac14 - \eps) n \ge \tfrac18 n$. As before, we can find a set $D_x$ of at least $\xi n$ distinct differences of $4$-APs contained in $S^-_A(x)$.
Hence, applying Lemma~\ref{lem:chebyshev} with $\lambda = \ell$ and $\kappa = \ceil{\xi n}$, we find that $p = \omega( (\ell n)^{-1/2} )= \omega( (tn)^{-1/2} )$ suffices to ensure $A \cup [n]_p$ is whp Schur. As $p(n,t)=t^{-1}\geq (tn)^{-1/2}$, this completes the proof in this case. 

\paragraph{Case II.2: $A'$ consists of large integers}

In this case, $\min A' > \card{A'}$. Since $\card{A'} \geq ( \tfrac12 - \eps )n$, it follows that $\min A' > (\tfrac 12 - \eps) n$. Thus, if we write $M = \left[\ceil*{\tfrac{n}{2}}, n \right] \setminus A$, we have
\[ \card{M} = \card*{ \left[ \ceil*{\tfrac{n}{2}}, n \right] \setminus A } \le \card*{ \left[ \ceil*{( \tfrac12 - \eps) n}, n \right] \setminus A'} = \card*{ \left[ \ceil*{(\tfrac12 - \eps)n}, n \right]} - \card{A'} < 2 \eps n. \]
Let $m := \card{M}$. Since $\card{A} = \ceil{ \tfrac{n}{2} } + t$, we must have $\card*{A \cap \left[1, \floor*{\tfrac{n}{2}} \right]} = m + t$.
Let $X$ be the set consisting of the $\floor*{\tfrac13 (m+t)}$ smallest elements in $A$, and consider $x \in X$. Observe that \[S^+_A(x) = \{ y \in A : x + y \in A \} = A \cap (A - x),\] where $A - x = \{ z - x : z \in A \}$.

If $x \le \tfrac{n}{2} - 3(m+t)$, consider the interval $I_1 := \left[ \ceil*{\tfrac{n}{2}}, \ceil*{\tfrac{n}{2}} + 3(m+t) \right]$. We then have \[S^+_A(x) \cap I_1 = A \cap (A-x) \cap I_1 = I_1 \setminus \left( M \cup (M-x) \right),\] and so at most $2m$ elements of $I_1$ can be missing from $S^+_A(x)$. Therefore $S^+_A(x)$ contains at least $m+3t$ elements out of an interval of length $3(m+t)+1$, and hence by Theorem~\ref{thm:APsupersat} (which we may apply over any interval, as arithmetic progressions are translation-invariant), there is some $\xi = \xi(\tfrac13) > 0$ and a set $D_x$ of size $3 \xi (m+t)$, such that for each $d \in D_x$ there is a $4$-AP in $S^+_A(x)$ with common difference $d$.

Otherwise, we must have $\tfrac{n}{2} - 3(m+t) < x \le \tfrac{n}{2} - \tfrac23 (m+t)$ (since there are $\ceil*{\tfrac23(m+t)}$ elements of $A \setminus X$ that are at most $\floor{\tfrac{n}{2}}$, $x$ cannot be any larger). This time, consider the interval $I_2 := \left[ x, \ceil*{\tfrac{n}{2}} + \ceil*{ \tfrac23(m+t)} \right]$. We have
\begin{align*}
I_2 \setminus S^+_A(x) = I_2 \setminus \left(A \cap (A - x) \right) &= \left( I_2 \setminus A \right) \cup \left( I_2 \setminus (A - x) \right) \\ &= \left( I_2 \setminus A \right) \cup \left( \left( (I_2 + x) \setminus A \right) - x \right).     
\end{align*} 
In this case, by our choice of $\eps$ and since $x$ is large and the interval $I_2$ is small, no missing element of $A$ can contribute to both $I_2 \setminus A$ and $I_2 \setminus (A - x)$ simultaneously. We therefore have
\begin{align*}
&	\card*{I_2 \setminus S^+_A(x)} = \card*{\left(I_2 \cup (I_2 + x) \right) \setminus A} \\
	&\le \card*{ \left[ x, \floor*{\tfrac{n}{2}} \right] \setminus A} + \card*{ \left( \left[ \ceil*{\tfrac{n}{2}}, \ceil*{\tfrac{n}{2}} + \ceil*{\tfrac23(m+t)} \right] \cup \left[ 2x, \ceil*{\tfrac{n}{2}} + \tfrac23(m+t) + x \right] \right) \setminus A} \\
	&\le \card*{ \left[ x, \floor*{\tfrac{n}{2}} \right] \setminus A} + \card*{ \left[ \ceil*{\tfrac{n}{2}}, n \right] \setminus A} \\
	&= \card*{ \left[ x, \floor*{\tfrac{n}{2}} \right] \setminus A}  + m \\
	&\le \tfrac{n}{2} - (x-1) - \tfrac23(m+t)  + m < \tfrac{n}{2} - x + \tfrac13 (m + t),
\end{align*}
where the penultimate inequality follows from the fact that there are at least $\ceil*{\tfrac23(m+t)}$ elements of $A$ that are larger than $x$ and at most $\floor*{\tfrac{n}{2}}$. Therefore 
\[ \card*{S^+_A(x) \cap I_2} > \card{I_2} - \left( \tfrac{n}{2} - x + \tfrac13(m+t) \right) \ge \tfrac13 (m+t). \]
Since $\card{I_2} \le 4 (m+t)$, it follows that $S^+_A(x)$ is dense in $I_2$, and so we may again apply Theorem~\ref{thm:APsupersat} to find some $\xi = \xi_{\ref{thm:APsupersat}}(\tfrac{1}{12}) > 0$ and a set $D_x$ of size at least $\xi \card{I_2} \ge \tfrac43 \xi (m+t)$ such that, for every $d \in D_x$, there is a $4$-AP in $S^+_A(x)$ with common difference $d$. 

We may therefore apply Lemma~\ref{lem:chebyshev} with $\lambda = \floor*{\tfrac13(m+t)}$ and $\kappa = \ceil*{\tfrac43 \xi (m+t)}$ to deduce that having $p = \omega \left( (m+t)^{-1} \right)$ ensures $A \cup [n]_p$ is whp Schur. Thus, our choice of $p = \omega \left( t^{-1} \right)$ is sufficient.
\end{proof}

\vspace{1em}

We now complete the proof of the $1$-statement of Theorem~\ref{thm:main_dense} by proving Lemma~\ref{lem:chebyshev}.

\begin{proof}[Proof of Lemma~\ref{lem:chebyshev}]
Given some $x \in X$ and $d \in D_x$, suppose first that $d$ is the common difference of a $4$-AP in $S^+_A(x)$. Let $a$ be the first term of such a $4$-AP. Then $\{a, a+d, a+2d, a+3d\} \subseteq S^+_A(x) \subseteq A$, and thus, by definition of $S^+_A(x)$, we also have \[\{ a+x, a+x+d, a+x + 2d, a+x+3d\} \subseteq A.\] In this case, we define $P(x,d) = \{d, x+d\}$. Note that $P(x,d) \subseteq [n]$, since $1 \le d \le x + d \le a + x + d \le n$, where the final inequality follows from the fact that $a + x + d \in A \subseteq [n]$. We further have $L_1(a,x,d) \setminus A \subseteq P(x,d)$. Since, by Proposition~\ref{prop:11set}, $L_1(a,x,d)$ is Schur, it follows that $A \cup [n]_p$ will be Schur if $P(x,d) \subseteq [n]_p$.

If, instead, $d$ is the common difference of a $4$-AP in $S^-_A(x)$ and $x\neq 2d$, letting $a$ be the first term of such a $4$-AP, then $\{a, a+d, a+2d, a+3d\} \subseteq S^-_A(x) \subseteq A$ and, by definition of $S^-_A(x)$, \[\{x-a-3d, x-a-2d, x-a-d, x-a\} \subseteq A.\] Here, we define $P(x,d) = \{d, x-d\}$, which is again easily seen to be contained in $[n]$. Then we have $L_2(a,x,d) \setminus A \subseteq P(x,d)$, and so, since $L_2(a,x,d)$ is Schur, it once more suffices to have $P(x,d) \subseteq [n]_p$.

Note that the map $(x,d) \mapsto P(x,d)$ is at most three-to-one; given a pair $\{u,v\}$ in the image with $u \le v$, it either takes the form of $\{d, x + d \}$ with $d = u$ and $x = v-u$, or the form $\{d, x-d\}$ with $x = u + v$ and $d = u$ or $d = v$. Hence, since there are at least $\tfrac12 \lambda \kappa$ pairs $(x,d)$ with $x\neq 2d$ (the factor of $\tfrac12$ comes from ignoring the pairs with $x=2d$), there are at least $\tfrac16 \lambda \kappa$ distinct pairs $P(x,d)$ whose appearance in $[n]_p$ would make $A \cup [n]_p$ Schur. Moreover, as we ignored cases in which $x=2d$, all the pairs $P(x,d)$ are indeed sets of size 2. Let $Y$ be the random variable counting how many of these pairs are contained in $[n]_p$. We have 
\[ \mb{E}[Y] \geq  \tfrac16 \lambda \kappa p^2 = \omega(1), \]
since $p = \omega \left( (\lambda \kappa)^{-1/2} \right)$.

To bound the variance of $Y$, note that the events of $P(x,d)$ and $P(x',d')$ appearing in $[n]_p$ are independent unless $P(x,d) \cap P(x',d') \neq \emptyset$. Furthermore, a given element $u$ can be in at most $4 \lambda$ pairs $P(x,d)$; once we specify which type of pair it is, and which role $u$ plays in the pair, each pair determines a unique $x \in X$. Thus, fixing a pair $P(x,d)$, it follows that there are at most $8\lambda$ other pairs which intersect it. Any such intersecting pair of pairs consists of a total of three elements, and they all appear in $[n]_p$ with probability $p^3$. Thus 
\[ \Var (Y) \le 8 \lambda^2 \kappa p^3 = o(\mb{E}[Y]^2), \]
since $p = \omega(\kappa^{-1})$.

Hence it follows from Chebyshev's Inequality (Theorem~\ref{thm:chebyshev}) that $\mb{P}(Y = 0) = o(1)$. That is, $[n]_p$ will whp contain some pair $P(x,d)$, and thus $A \cup [n]_p$ will be Schur.
\end{proof}

\subsubsection{Proof of Theorem~\ref{thm:stability}}

Suppose $A \subseteq [n]$ is a set of size $\ceil*{\tfrac{n}{2}} + t$, and $q = \omega \left(n^{-1} \right)$ is such that $A \cup [n]_q$ is whp not Schur. From the proof of Theorem~\ref{thm:main_dense}, we have that if $A$ fell under Case I then $p=\omega(n^{-1})$ would be sufficient to ensure that $A \cup [n]_p$ is whp Schur. Hence, we must have that $A$ falls into Case II. Now if $A$ falls into Case II.1, the proof of Theorem~\ref{thm:main_dense} shows that if $A$ contains $\ell$ even numbers and $p=\omega((\ell n)^{-1/2})$ then whp $A\cup [n]_p$ is Schur. In this case we must thus have $q=O((\ell n)^{-1/2})$ and so $\ell=O(q^{-2}n^{-1})$.

On the other hand, if $A$ falls under Case II.2, the proof shows that, for $ m := \card*{\left[ \ceil*{\tfrac{n}{2}}, n\right] \setminus A}$, $p = \omega \left( (m+t)^{-1} \right)$ suffices to make $A \cup [n]_p$ Schur whp. Thus we must have $q = O \left((m+t)^{-1} \right)$, which in particular implies $m = O \left( q^{-1} \right)$, proving the stability result. \qed

\section{Sparse base sets}
\label{sec:sparse}

In this section we prove Theorem~\ref{thm:main_sparse}, starting by proving the $0$-statement in Section~\ref{sec:0sparse}. In Section~\ref{sec:containers} we then introduce containers for colourings, which will be a key tool in proving the $1$-statement of Theorem~\ref{thm:main_sparse}, which we carry out in Section~\ref{sec:1sparse}.

\subsection{The $0$-statement of Theorem~\ref{thm:main_sparse}}
\label{sec:0sparse}

As with the proof of the $0$-statement of Theorem~\ref{thm:main_dense}, we will take $A$ to be the set of the largest integers in $[n]$, which we denote by $A_s:=[n-s+1,n]$. However, in contrast to the setting of dense base sets, our proof here is non-constructive. We obtain a contradiction by assuming that the random perturbation of $A_s$ \emph{is} Schur and appealing to a minimal Schur subset to derive that  the random set must then contain one of a collection of substructures, all of which do not appear whp. This  approach to proving $0$-statements dates back to the pioneering work for finding thresholds for Ramsey properties of graphs \cite{luczak1992ramsey,rodl1993lower,rodl1995threshold} as well as the proof of the $0$-statement of Theorem~\ref{thm:random} due to Graham, R\"odl and Ruci\'nski~\cite{graham1996schur}.

Recall from Section~\ref{sec:Notation} the definition of $\mc H_{\textrm{Schur}}(X)$ for a set $X\subseteq [n]$ and note that the edges of $\mc H_{\textrm{Schur}}(X)$ have cardinality two (in the case of degenerate Schur triples) or three. Note further that finding a Schur colouring of a set $X \subseteq [n]$ is equivalent to finding a proper two-colouring of $\mc H_{\textrm{Schur}}(X)$; that is, a two-colouring without monochromatic edges. Therefore, to prove the $0$-statement of the theorem, we need to show that if $p = o((sn)^{-1/3})$, then whp, after we add the perturbation $P := [n-s]_p$ to the set $A_s$, we can properly two-colour $\mc H_{\textrm{Schur}}(A_s \cup P)$. We will in fact prove a bit more, showing that one can find such a colouring where all elements in $A_s$ are coloured blue.

Now suppose that $\mc H_{\textrm{Schur}}(A_s \cup P)$ does not admit a proper two-colouring where all elements of $A_s$ are coloured blue, and let $\mc H_{min}=\mc{H}_{min}(s,P) \subseteq \mc H_{\textrm{Schur}}(A_s\cup P)$ be an edge-minimal subgraph with this property. That is, $\mc H_{min}$ does not have such a two-colouring, but the hypergraph obtained by removing any edge $h$ from $\mc H_{min}$ does. We now explore the structure of $\mc H_{min}$.

\begin{lem} \label{lem:extension}
Suppose $P$ is such that $\mc H_{\textrm{Schur}}(A_s \cup P)$ does not admit a proper two-colouring where all elements of $A_s$ are coloured blue. Then the minimal subgraph $\mc H_{min} = \mc{H}_{min}(s,P)$ defined above has the following property. 
Every edge $h \in \mc H_{min}$ contains at least one element from $P$, and, for every $x \in h \cap P$, there is an edge $h' \in \mc H_{min}$ such that $h \cap h' = \{x \}$. Moreover, if $h$ contains an element from $A_s$, then there exists such an edge $h'$ that does not contain elements from $A_s$.
\end{lem}

\begin{proof}
First observe that, as $s \le \floor{\tfrac{n}{2}}$, the set $A_s$ is sum-free. Hence, $A_s$ is an independent set in $\mc H_{\textrm{Schur}}(A_s \cup P)$, and every edge of $\mc H_{min} \subseteq \mc H_{\textrm{Schur}}(A_s \cup P)$ must contain an element of $P$.

Now, by minimality of $\mc H_{min}$, we know that there is a proper colouring of $\mc H_{min} \setminus \{h\}$ where all elements in $A_s$ are blue. As $\mc H_{min}$ itself does not admit such a colouring, $h$ must be monochromatic under this colouring. If we swap the colour of $x$, then $h$ will no longer be monochromatic, so we must create another monochromatic edge, say $h'$. As $x$ was the only element to change colour, we must have $x \in h'$, and as $x$ did change colour, $h'$ must be a different colour than $h$ was. Hence $h$ and $h'$ cannot have any other elements in common, and $h \cap h' = \{x\}$.

To establish the final assertion, observe that if $h \cap A_s$ is non-empty, $h$ must have been coloured blue. Thus the edge $h'$ must be coloured red after recolouring $x$, and hence cannot contain any element from $A_s$.
\end{proof}

While the above result holds for any outcome of the random set $P$, our next proposition gives some additional structure that holds whp.

\begin{prop} \label{prop:minhypergraph}
The random set $P$ is such that the following holds whp. If $P$ is such that $\mc H_{\textrm{Schur}}(A_s \cup P)$ does not admit a proper two-colouring where all elements of $A_s$ are coloured blue, then the hypergraph $\mc H_{min}=\mc H_{min}(s,P)$ has the following properties.
\begin{enumerate}
    \item[(a)] The hypergraph $\mc H_{min}$ is three-uniform.
    \item[(b)] Every edge of $\mc H_{min}$ contains at most one element from $A_s$.
    \item[(c)] The hypergraph $\mc H_{min}$ is linear. That is, any pair of distinct edges of $\mc H_{min}$ intersect in at most one vertex. 
\end{enumerate}
\end{prop}

\begin{proof}
We start with part (a). As noted earlier, all edges of $\mc H_{\textrm{Schur}}(A_s \cup P)$ have cardinality either two or three, with the two-edges $h$ of the form $\{x,2x\}$. If $2x \le n-s$, then $h \subseteq P$, and so each of its elements appears independently with probability $p$. Thus, the expected number of such edges in $\mc H_{\textrm{Schur}}$ is bounded by $np^2$. Since $p=o( n^{-1/2})$, this is $o(1)$ and so whp none of these edges appear.

On the other hand, if $2x \ge n-s+1$, then $2x \in A_s$, and there are at most $s$ choices for $2x$. However, $x \le \tfrac{n}{2} < n-s+1$, and so $x \in P$. By Lemma~\ref{lem:extension}, there is an edge $h'$ that meets $h$ only in $x$, and this edge must be fully contained in $P$. By the preceding paragraph, $\card{h'} = 3$. Since the element $x$ participates in at most $n$ Schur triples, there are at most $n$ choices for $h'$, and the probability that $h' \subseteq P$ is $p^3$. Thus, the probability that any such $h'$ exists in $\mc H_{min}$ (and, therefore, that $h$ does) is at most $np^3$. Taking a union bound over the possible choices for $h$, we find the probability of the existence of such a two-edge is at most $snp^3 = o(1)$.

We now turn to part (b). We know from Lemma~\ref{lem:extension} that there is no edge fully contained in $A_s$, so we need only show that there are no edges in $\mc H_{min}$ containing two elements of $A_s$. Suppose $h = \{x,y,z\}$ were such an edge, with $x + y = z$, where necessarily $x \in P$ and $y,z \in A_s$. We must then have $x \in [s-1]$, and so there are fewer than $s$ choices for $x$. By Lemma~\ref{lem:extension}, there is an edge $h' \subseteq P$ with $h \cap h' = \{x\}$. For each choice of $x$, there are at most $n$ choices for $h'$, each of which appears with probability $p^3$. Hence, taking a union bound over all choices of $x$ and $h'$, the probability that such a supporting edge exists in $\mc H_{min}$ is at most $snp^3 = o(1)$.

For part (c), suppose we have two edges $h = \{x,y,z\}$ and $h' =  \{x,y,w\}$. Note that by part (b), it suffices to focus on the case that both $h$ and $h'$ contain at most one element of $A_s$. Indeed we have that  whp there are no pairs of edges $h,h'$ with either $h$ or $h'$  containing more than one element from $A_s$. Now we first consider the case $h \cup h' \subseteq P$. There are at most $n$ choices for $x$ and at most $n$ choices for $y$, but once these are fixed, there are only constantly many choices for $z$ and $w$. As each element appears in $P$ independently with probability $p$, the probability of finding such a configuration is at most $O(n^2p^4) = o(1)$.

Next we handle the case $x \in A_s$ (the case $y \in A_s$ is symmetric). There are at most $s$ choices for $x$ and at most $n$ choices for $y$. Once this pair is fixed, there are again only a constant number of choices for $z$ and $w$. Moreover, as $h$ and $h'$ have at most one vertex in $A_s$, we have that $y, z, w \in P$. Hence, taking a union bound over all possible such configurations, the probability that one appears is at most $snp^3 = o(1)$.

Finally, suppose we have $z \in A_s$. In this case, we must have $x,y \in P$, with the relation $x + y = z$. We must further have $w = \card{y-x}$, and so $w \in P$ as well. There are then at most $s$ choices for $z$ and $n$ choices for $x$, after which $y$ and $w$ are determined. The probability of finding such a configuration is therefore at most $snp^3 = o(1)$, completing the proof of (c).
\end{proof}

\begin{figure}[b]
    \centering
    \includegraphics[scale = 0.8]{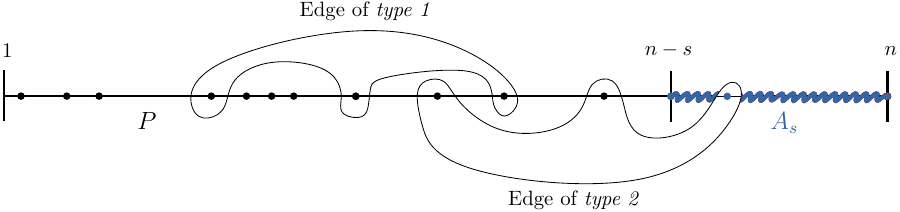}
    \caption{An edge of type 1 and an edge of type 2}
    \label{fig:edge_type}
\end{figure}

In the following it will help to differentiate between two kinds of edges that appear in $\mc H_{min}$. We call an edge \emph{type 1} if it is fully contained in $P$, and \emph{type 2} if it contains exactly one vertex from $A_s$. See Figure~\ref{fig:edge_type} for an example. Also, in what follows, a \emph{loose cycle} in a $3$-uniform hypergraph $\mc H$ is a collection of $\ell\geq 3$ edges $e_1,\ldots,e_\ell\in E(\mc H)$ such that for all $i\in [\ell]$, $|e_i\cap e_j|=1$ for $j=i-1,i+1 \mod \ell$, and $e_i\cap e_j= \emptyset$ if $j\neq i-1,i,i+1 \mod \ell$.  

\begin{lem}
The random set $P$ is such that the following holds whp. If $P$ is such that $\mc H_{\textrm{Schur}}(A_s \cup P)$ does not admit a proper two-colouring where all elements of $A_s$ are coloured blue, then there is a loose cycle in $\mc H_{min}=\mc H_{min}(s,P)$ that has at most one pair of consecutive type 2 edges and such that all degree 2 vertices in the cycle belong to $P$.
\label{lem:cycle}
\end{lem}

\begin{proof}
Using Lemma~\ref{lem:extension} and Proposition~\ref{prop:minhypergraph} we construct a walk in $\mc H_{min}$. This walk must eventually repeat a vertex, at which point it creates a cycle. 

To build the walk, start from an arbitrary edge $h_0 \in \mc H_{min}$. By Lemma~\ref{lem:extension}, there must be some $x_0 \in h_0 \cap P$. Now suppose we have already taken $i \ge 0$ steps in the walk, and have some edge $h_i \in \mc H_{min}$ and vertex $x_i \in h_i \cap P$. By Lemma~\ref{lem:extension}, there is some edge $h_{i+1} \in \mc H_{min}$ such that $h_i \cap h_{i+1} = \{x_i\}$. Furthermore, by Proposition~\ref{prop:minhypergraph}(a), $\card{h_{i+1}} = 3$, and by Proposition~\ref{prop:minhypergraph}(b), we can find some $x_{i+1} \in h_{i+1} \cap P$ that is distinct from $x_i$. Hence, we can extend the walk with the edge $h_{i+1}$, and proceed to the next step using the vertex $x_{i+1}$. 

We repeat this process until the last edge $h_t$ added contains a vertex, apart from $x_{t-1}$, that we have seen previously. That is, there is some $r \le t-1$ and some $y \in h_t \setminus \{ x_{t-1} \}$ such that $y \in h_t \cap h_r$. If there are multiple choices for the index $r$, we choose the larger one. Note that, by Proposition~\ref{prop:minhypergraph}(c), $\mc H_{min}$ is linear, and so we must in fact have $r \le t-2$, as $h_{t-1}$ and $h_t$ already share the vertex $x_{t-1} \neq y$. This means the edges $h_r, h_{r+1}, \hdots, h_{t-1}, h_t$ form a loose cycle in $\mc H_{min}$.

Lemma~\ref{lem:extension} guarantees that if an edge $h_i$ is of type 2, then the following edge $h_{i+1}$ will be of type 1. Hence, we do not have two consecutive type 2 edges in the walk. Note that it is possible that the first and last edges $h_r$ and $h_t$ might both be of type 2. Thus, in the cyclic ordering of the edges, there is at most one pair of consecutive type 2 edges.
\end{proof}

Using these properties, we can prove the $0$-statement of Theorem~\ref{thm:main_sparse}.

\begin{proof}[Proof of the $0$-statement of Theorem~\ref{thm:main_sparse}]
Assume that $P$ satisfies the conclusion of Lemma~\ref{lem:cycle} (which holds whp) and suppose for a contradiction that $P$ is such that $\mc H_{\textrm{Schur}}(A_s \cup P)$ does not admit a proper two-colouring where all elements of $A_s$ are coloured blue. Then, appealing to Lemma~\ref{lem:cycle}, we can construct a loose cycle in the associated minimal hypergraph $\mc H_{min}=\mc H_{min}(s,P) \subseteq \mc H_{\textrm{Schur}}(A_s \cup P)$ that has at most one pair of consecutive type 2 edges. In the following we show that whp such a cycle does not exist in $\mc H_{min}$, and hence there must exist a Schur colouring of $A_s\cup P$.

We start by ruling out cycles without any consecutive type 2 edges. If the cycle has length $\ell$, for some $\ell \ge 3$, we label the edges in cyclic order as $e_1, e_2, \hdots, e_\ell$, choosing an ordering where $e_{\ell}$ is of type 1. We can then label the vertices of the cycle in such a way that, for each $1 \le i \le \ell-1$, $e_i = \{x_{i-1}, y_i, x_i\}$, while $e_\ell = \{x_{\ell-1}, y_{\ell}, x_0\}$. Note that, by Lemma \ref{lem:cycle}, each degree-two vertex $x_i$ belongs to $P$, while $y_i \in A_s$ if and only if the edge $e_i$ is of type 2.
We can now bound the expected number of such cycles. We start by choosing the vertex $x_0$, for which there are at most $n$ choices, each appearing with probability $p$. Hence, this choice contributes a factor of $np$.

Now, consider the choices for the edges $e_i$, where $1 \le i \le \ell-\varepsilon$, where $\varepsilon = 1$ if $e_{\ell-1}$ is of type 1 and $\varepsilon = 2$ if $e_{\ell-1}$ is of type 2. If $e_i$ is of type 1, there are at most $n$ choices of $y_i$, and then at most two choices for $x_i$ that complete a sum with $x_{i-1}$ and $y_i$. The elements $y_{i}$ and $x_i$ appear in $P$ independently, each with probability $p$. Thus a type 1 edge contributes a factor of at most $2np^2$.

On the other hand, if $e_i$ is of type 2, then there are at most $s$ choices of $y_i \in A_s$, and then $x_i = y_i - x_{i-1} \in P$ is determined, and appears in $P$ with probability $p$. Thus, this edge contributes a factor of $sp$. However, we then know the next edge $e_{i+1}$ is of type 1, and contributes a factor of $2np^2$. We group these two factors: every type 2 edge, together with its subsequent type 1 edge, contributes a factor of $2snp^3$.

Hence, in total, each such intermediate extension contributes a factor of at most $(2np^2 + 2snp^3)$.
This leaves us with the task of closing the cycle. If $e_{\ell-1}$ is of type 1, then we have already accounted for it previously, and need only consider the choice of $e_\ell$. Recall that we chose our labelling so that $e_{\ell}$ is of type 1. The vertices $x_{\ell-1}$ and $x_0$ are already fixed, and there are at most two choices for the final vertex $y_\ell \in P$, which appears with probability $p$. Thus, we gain a factor of $2p$ in this case.
If, on the other hand, $e_{\ell -1}$ is of type 2, then the same arguments as before show that we gain a factor of at most $sp$ for the edge $e_{\ell-1}$, and a factor of $2p$ for the final type 1 edge $e_{\ell}$. Hence, in this case, we gain a factor of $2sp^2$.

In total, closing the cycle contributes a factor of at most $2p(1 + sp)$.
When we sum over all possible cycles, therefore, we can bound the expectation by
\[ np \cdot \sum_{k \ge 0} \left(2np^2 + 2snp^3\right)^k \cdot 2p(1 + sp) = \sum_{k \ge 0} \left( 2np^2 + 2snp^3 \right)^{k+1}, \]
where $k$ represents the number of intermediate extensions in forming the cycle. By virtue of the fact that $p=o( (sn)^{-1/3})$ and $p=o( n^{-1/2})$, this sum is $o(1)$, and so whp we do not have any such cycles.

\medskip

This leaves us with those cycles whose initial and final edges $h_r$ and $h_t$ (adopting the notation from the proof of Lemma~\ref{lem:cycle}) are both of type 2. We further split into two subcases, based on whether their common vertex lies in $A_s$ or $P$.
In the former setting, we index the edges in cyclic order starting with $h_r$ and closing the cycle with $h_t$, so that we have $h_r = e_1, e_2, e_3, \hdots, e_{\ell - 1}, e_{\ell} = h_t$. We label the vertices within the edges similarly to before, except when it comes to $e_{\ell}$, as the common vertex with $e_1$ will be the $A_s$ vertex $y_1$. That is, $e_i = \{x_{i-1}, y_i, x_i\}$ for all $i \in [\ell - 1]$, while $e_{\ell} = \{x_{\ell - 1}, y_1, x_{\ell} \}$. As before, each vertex $x_i$ lies in $P$, while $y_i$ lies in $P$ if $e_i$ is of type 1, and in $A_s$ otherwise.

We can then bound the expected number of such cycles just as we did before. For the initial edge $e_1$, we have at most $n$ choices for $x_0 \in P$, at most $s$ choices for $y_1 \in A_s$, and then $x_1 \in P$ is determined uniquely. The vertices $x_0$ and $x_1$ appear in $P$ independently, each with probability $p$, and hence the contribution of $e_1$ to the expectation is at most a factor of $snp^2$.
Now, since the final edge $e_{\ell}$ is of type 2, its predecessor $e_{\ell-1}$ must be of type 1, and hence the intermediate extensions account for the edges $e_2$ through to $e_{\ell - 1}$. As before, each extension contributes a factor of at most $(2np^2 + 2snp^3)$.

Finally, when closing the cycle with the edge $e_{\ell}$, we have already fixed the elements $x_{\ell-1}$ and $y_1$, and so the edge $e_{\ell}$ is determined. However, $x_{\ell} \in P$ is a new vertex, and appears with probability $p$. We therefore collect a factor $p$ for the final edge $e_{\ell}$.
Putting this all together, the expected number of cycles of this type can be bounded by
\[ snp^2 \cdot \sum_{k \ge 0} \left( 2np^2 + 2snp^3 \right)^k \cdot p = snp^3 \sum_{k \ge 0} \left( 2np^2 + 2snp^3 \right)^k = o(1), \]
and so whp we do not have any such cycles.

\medskip

For the latter subcase, where the edges $h_r$ and $h_t$ meet in a vertex $x_1 \in P$, we shall instead have these two edges be the first two in our cyclic ordering. We label the edges of the cycle as $e_1 = h_t, e_2 = h_r, e_3, e_4, \hdots, e_{\ell - 1}, e_{\ell}$, and we label the vertices within the edges as before, with $e_i = \{ x_{i-1}, y_i, x_i \}$ for $i \in [\ell - 1]$, and $e_{\ell} = \{x_{\ell-1}, y_{\ell}, x_0\}$. Note that $e_{\ell}$ must be a type 1 edge, as we cannot have another pair of consecutive type 2 edges.

Furthermore, observe that the vertex $x_1$ is only used in the cycle as the intersection between $e_1$ and $e_2$, two edges of type 2. By Lemma~\ref{lem:extension}, we are guaranteed the existence of another edge, say $f$, of type 1, such that $f \cap e_1 = \{x_1\}$. If $f$ were to contain another vertex from the cycle, that would create a cycle without a pair of consecutive type 2 edges, but we previously showed that such cycles do not exist. Hence, the vertices in $f \setminus \{x_1\}$ must be new.

We now bound the expected number of copies of these cycles, together with the pendant edge $f$. There are $n$ choices for the vertex $x_1$, which appears with probability $p$. There are then $s$ choices for each edge $e_1$ and $e_2$, and their other $P$-vertices, $x_0$ and $x_2$, appear independently with probability $p$ each. Finally, there are at most $n$ choices for the edge $f$, and the two vertices in $f \setminus \{x_1\}$ also each appear with probability $p$. Thus, the initial constellation of edges $e_1, e_2$ and $f$ contributes a factor of at most $s^2 n^2 p^5$ to the expectation.

Since $e_2$ is of type 2, the edge $e_3$ must be of type 1. Hence, for the edges $e_3, e_4, \hdots, e_{\ell - 1}$, every type 2 edge is preceded by a type 1 edge, and so we can again\footnote{Previously we paired a type 2 edge with the type 1 edge succeeding it, but the calculations here are identical.} bundle these together when bounding the expectation. Then, as before, each intermediate extension provides a factor of at most $(2np^2 + 2snp^3)$, while closing the cycle with the type 1 edge $e_{\ell}$ gives an additional factor of $2p$. Thus the expectation can be bounded by
\[ s^2n^2p^5 \cdot \sum_{k \ge 0} \left( 2np^2 + 2snp^3 \right)^k \cdot 2p = 2 \left( snp^3 \right)^2 \sum_{k \ge 0} \left( 2np^2 + 2snp^3 \right)^k = o(1), \]
and so again we do not have any such cycle whp.

In summary, we find that whp $\mc H_{min}$ does not contain any cycle that would result from Lemma~\ref{lem:cycle}, which means that our initial assumption that $A_s \cup [n]_p$ does not have a Schur colouring with all elements of $A_s$ coloured blue, cannot hold. This completes the proof.
\end{proof}

\subsection{Containers for colourings}~\label{sec:containers}
For the $1$-statement of Theorem~\ref{thm:main_sparse}, we fix a set $A \subseteq [n]$ of the appropriate size, and wish to show that when $p$ is sufficiently large, then whp a random perturbation $P \sim [n]_p$ is such that $A \cup P$ is Schur. Roughly speaking, the idea is that, for a given colouring of $[n]$, the perturbation $P$ is very likely to contain elements that, in combination with $A$, form a monochromatic Schur triple. Unfortunately, there are far too many potential colourings of $[n]$, rendering the union bound ineffective.

To resolve this issue, we make use of hypergraph containers, introduced by Saxton and Thomason~\cite{saxton2015hypergraph} and Balogh, Morris and Samotij~\cite{BMS15}, which have been successfully applied to a wide range of problems in combinatorics. In our setting, the key idea is to group similar colourings into so-called containers, and then show that the random set is unlikely to fit with not just a given colouring, but also the container at large. As the number of containers will be much smaller, we will then be able to proceed with a union bound and obtain the desired result.

\medskip

To put things on a formal footing, we define a colouring hypergraph $\mc H_A$ that will encode the colourings of subsets of $[n]$. The vertex set $V(\mc H_A)$ is the disjoint union of two copies of $[n]$, which we call $V_R$ and $V_B$. Colouring the element $i \in [n]$ red will then be represented by the vertex $i_R \in V_R$, while colouring it blue will be represented by $i_B \in V_B$. Thus, given a subset $S \subseteq [n]$ and a colouring $\varphi: S \to \{ \text{red}, \text{blue} \}$, we can identify $\varphi$ with the vertex set $\{i_R \in V_R: i \in S, \varphi(i) = \text{red} \} \cup \{i_B \in V_B: i \in S, \varphi(i) = \text{blue} \}$.

The edges of the hypergraph will correspond to coloured configurations that cannot appear in a Schur colouring of $A \cup P$. Indeed, given some $a \in A$, let $x,y,w,z \in [n]$ be such that the sets $\{a,x,y\}$ and $\{a,z,w\}$ both host Schur triples. If we were to colour $x$ and $y$ red and colour $z$ and $w$ blue, then assigning either colour to $a$ creates a monochromatic sum. Hence, this colouring of the four elements $x,y,z$ and $w$ can be forbidden, motivating our definition of the edge set of $\mc H_A$: 
\begin{align*}
 E(\mc H_A) = \Big\{ \{ x_R, y_R,  & z_B,  w_B\}    :   \exists  \; a \in A    \text{ for which } \\  &  \{a,x,y\}  \text{ and } \{a,z,w\} 
\text{ host non-degenerate Schur triples} \Big\}.     
\end{align*}
We restrict here only to non-degenerate Schur triples to ease the analysis. Given an edge $e=\{x_R, y_R, z_B, w_B\}$, we call the associated element $a \in A$ the \emph{target}\footnote{When $\{x,y\} = \{z,w\}$, the target of the edge need not be unique --- it could be either the difference or the sum of the pair. In such a case, when referring to the target of the edge, we arbitrarily choose one such target.}
 of the edge $e$.
 
It thus follows that if $A \cup P$ admits a Schur colouring $\varphi$, then $\varphi \subseteq V(\mc H_A)$ must be an independent set. The theory of hypergraph containers asserts that when the edges of a hypergraph are well-distributed, in the sense that no set of vertices is contained in a disproportionally large number of edges, then all its independent sets can be grouped together into a small number of containers, each of which induces few edges.

To make the condition on the hypergraph precise, we must define the codegree function. Given an $r$-uniform hypergraph $\mc H$ on $N$ vertices and with average (vertex) degree $d$, some vertex $v \in V(\mc H)$, and some $1\leq j \le r$, we define $d_j(v) = \max \left\{ d(\sigma) : v \in \sigma \subseteq V(\mc H), \card{\sigma} = j \right\}$, where $d(\sigma)$ is the number of edges containing $\sigma$. Then, given any $\tau > 0$, we set $\delta_j = \frac{\sum_v d_j(v)}{\tau^{j-1} N d}$, and define the co-degree function \[\delta(\mc H, \tau) = 2^{\binom{r}{2} - 1} \sum_{j=2}^r 2^{-\binom{j-1}{2}} \delta_j.\] With this notation in place, we can state a version of the hypergraph containers theorem due to Saxton and Thomason.

\begin{thm}[Container theorem, Corollary 3.6 in \cite{saxton2015hypergraph}]
\label{thm:container}
Let $\mc H$ be an $r$-uniform hypergraph on the vertex set $[N]$, and suppose that $0<\tau,\eps<1/2$ satisfy $\delta(\mc H,\tau)\le \eps/12r!$. Then there are constants $c=c(r)$ and $z \le c \log (1 / \eps)$ and a function $ \Psi \colon\mathcal{P}([N])^z \to \mathcal{P}([N])$ with the following properties. Let $\mathcal{T} =\{ (T_1,\ldots, T_z)\in \mathcal{P}([N])^z\colon |T_i|\le c\tau N, 1\le i\le z\}$, and let $\mathcal{C}=\{ \Psi(T)\colon T\in \mathcal{T}\}$. Then
\begin{enumerate}
    \item For every independent set $I$ there exists $T=(T_1,\ldots,T_z)\in \mathcal{T}$ with $I\subseteq \Psi(T)\in \mathcal{C}$,
    \item $e(\mc H[C])\le \eps e(\mc H)$ for all $C\in \mathcal{C}$, and
    \item $\log |\mc C|\le c\log (1 / \eps) \tau N \log (1/\tau)$.
\end{enumerate}
\end{thm}

Applying this to our hypergraph $\mc H_A$, we arrive at the following.

\begin{cor}
\label{cor:containers}
For every fixed $\eps > 0$ there is a constant $c = c_{\eps}$ such that, if $A \subseteq [n]$ is a set of size $s = \Omega \left( n^{1/2} \right)$, then there is a collection $\mc C$ of subsets of $V(\mc H_A)$ for which the following hold:
\begin{enumerate}
    \item For every $P \subseteq [n]$ and Schur colouring $\varphi$ of $A \cup P$, there is some $C \in \mc C$ such that $\varphi \subseteq C$.
    \item For every $C \in \mc C$, $e(\mc H_A[C]) \le \eps sn^2$.
    \item $\log \; \card{\mc C} \le c s^{-1/3} n^{2/3} \log n$.
\end{enumerate}
\end{cor}

\begin{proof}
In order to derive this from Theorem~\ref{thm:container}, we need to compute the codegree function of the hypergraph $\mc H_A$. To start, we count the edges of $\mc H_A$. Observe that for each of the $s$ choices of the target of the edge, there are at least $\tfrac12 n-1$ and at most $n$ choices for the red pair forming a Schur triple, with the same bounds holding for the blue pair. Moreover, every edge has at most 2 targets. Thus, we have $\tfrac1{16} sn^2 \le e(\mc H_A) \le sn^2$. As there are $2n$ vertices, it follows that the average degree $d$ satisfies $d \ge \tfrac18 sn$.

We next need to bound the quantities $d_j(v)$, $2 \le j \le 4$, from above. To this end, we define $\Delta_j$ to be the maximum degree of a set $\sigma$ of $j$ vertices, noting that $\Delta_j = \max \{d_j(v) : v \in V(\mc H_A)\}$.
When $j = 2$, there are two cases to consider. If the two vertices in $\sigma$ are from the same colour, then there can be at most two choices for the target of the edge. This in turn leaves at most $n$ choices for the pair of vertices of the opposite colour. Hence, for any $\sigma$ of this form, we have $d(\sigma) \le 2n$.
On the other hand, if $\sigma$ has one vertex of each colour, then there are at most $s$ choices for the target of the edge. Once the target is chosen, there are again at most two choices for the remaining vertex of each colour, and thus $d(\sigma) \le 4s \le 4n$. Hence we deduce $\Delta_2 \le 4n$.

When $j = 3$, observe that $\sigma$ must contain both vertices from one of the colours, and thus there are only at most two possibilities for the target of the edge. Once this is fixed, and since we already have one vertex from the other colour, there are again at most two choices for the missing vertex, and thus $\Delta_3 \le 4$. Finally, we trivially have $\Delta_4 = 1$, since each edge consists of four vertices.

We then have $\delta_j = \frac{\sum_v d_j(v)}{\tau^{j-1} 2n d} \le \frac{2n \Delta_j}{\tau^{j-1} 2n d} \le \frac{8\Delta_j}{\tau^{j-1} sn}$, and so 
\[ \delta(\mc H_A, \tau) = 32 \delta_2 + 16 \delta_3 + 4 \delta_4 \le \frac{2^{10}}{\tau s} + \frac{2^{9}}{\tau^2 sn} + \frac{2^5}{\tau^3 sn}. \]
Therefore there is a constant $c'$ depending on $\eps$ such that, if \[\tau \ge c' \max \{ s^{-1}, (sn)^{-1/2}, (sn)^{-1/3}, \}\] we have $\delta(\mc H_A, \tau) \le \tfrac{1}{288} \eps$. Since $s = \Omega \left( n^{1/2} \right)$, the bound simplifies to $\tau \ge c' (sn)^{-1/3}$.

We can then apply Theorem~\ref{thm:container} with this choice of $\tau$, and the corollary follows immediately. Indeed, the first conclusion follows from our observation that such a Schur colouring $\varphi$ corresponds to an independent set in $\mc H_A$, and hence must be contained in a container. The second conclusion is a consequence of our earlier calculation showing $e(\mc H_A) \le sn^2$. Finally, the bound on the number of containers comes from making the substitutions $\tau = c' (sn)^{-1/3}$ and $N = 2n$ in the corresponding bound from the theorem.
\end{proof}

\subsection{The $1$-statement of Theorem~\ref{thm:main_sparse}}
\label{sec:1sparse}

We will now use the containers of the previous section to prove the $1$-statement for sparse base sets. As previously discussed, the value of containers lies in the fact that when considering potential colourings, we can work on the level of the containers, which allows for a much more efficient union bound.

To make things precise, suppose we have a set of  containers $\mc C$ provided by Corollary~\ref{cor:containers}. We say that a set $P \subseteq [n]$ is \emph{compatible} with a container $C \in \mc C$, denoted $P \triangleleft C$, if there is some Schur colouring $\varphi$ of $A \cup P$ such that $\varphi \subseteq C$; that is, $\{ i_R: i \in A \cup P, \; \varphi(i) = \textrm{red} \} \cup \{ i_B : i \in A \cup P, \; \varphi(i) = \textrm{blue} \} \subseteq C$.
Note that if $P $ is such that $A \cup P$ is not Schur, then $A \cup P$ admits a Schur colouring, and thus by the first part of Corollary~\ref{cor:containers}, $P$ must be compatible with some container $C$. The key proposition below shows that this is exceedingly unlikely when $P\sim [n]_p$ and $p$ is sufficiently large.

\begin{prop}
There exists $\eps >0$ such that the following holds. Let $n$ and $s=s(n)$ be positive integers such that $s=\Omega(n^{1/2})$ and $s=o(n)$. Furthermore, let $A\subset [n]$ be a set of size $s$, let $\mc{H}_A$ be the colouring hypergraph as described in Section~\ref{sec:containers} and let $\mc{C}$ be the collection of containers corresponding to this hypergraph given by Corollary~\ref{cor:containers} with parameter $\eps$. Finally, let $p=p(n)$ such that $p=\omega((sn^{13})^{-1/27}\log n)$ and $p = O\left(n^{-1/2}\right)$ and let $P\sim [n]_p$. Then
\[\Pr\left[P\triangleleft C\right] = e^{-\omega\left(s^{-1/3}n^{2/3}\log n\right)},\] for any container $C\in \mc{C}$.
\label{prop:probability_container}
\end{prop}

Before giving the proof of this proposition, let us see how it implies the 1-statement of Theorem~\ref{thm:main_sparse}.

\begin{proof}[Proof of the 1-statement of Theorem~\ref{thm:main_sparse}]

Let $A$ be our base set with $|A|= s = \Omega(n^{1/2}) $ and let $p =\omega( (sn^{13})^{-1/27}\log n)$. If it is the case that $p=\omega\left(n^{-1/2}\right)$ then the conclusion follows directly from Theorem~\ref{thm:random} (as the random perturbation $P$ itself will be Schur whp) and so we may assume that $p=O\left(n^{-1/2}\right)$. Furthermore, if $s=\Omega(n)$, then the desired conclusion follows from Theorem~\ref{thm:posdense}, and so we may assume that $s=o(n)$. Now let $\mc H_A$ be the $4$-uniform colouring hypergraph described in Section~\ref{sec:containers} and let $P\sim [n]_p$ be the random perturbation. Let $\eps$ be the constant given by  Proposition~\ref{prop:probability_container} and apply the container lemma, or more precisely, Corollary~\ref{cor:containers}, with the parameter $\eps$.

Proposition~\ref{prop:probability_container} tells us that the probability that $P$ is compatible with a given container is at most $e^{-\omega\left(s^{-1/3}n^{2/3}\log n\right)}$. Moreover, if $A\cup P$ is not Schur, then there is a Schur colouring $\varphi$ of $A\cup P$ and, appealing to part 1 of Corollary~\ref{cor:containers}, $\varphi$ corresponds to a subset of one of the containers $C\in\mc C$ in $\mc H_A$.  Hence the event that $A\cup P$ is not Schur is contained in the event that $P$ is compatible with some container $C\in \mc C$. 
Taking a union bound over all containers we get that 
\[
\Pr[A\cup P \text{ is not Schur}] \le \sum_{C\in \mc C} \Pr[P\triangleleft C] \le
    e^{cs^{-1/3}n^{2/3}\log n }\: e^{-\omega\left(s^{-1/3}n^{2/3}\log n\right)} = o(1)\]
as required, where we used part $3$ of Corollary~\ref{cor:containers} to upper bound the number of containers. 
\end{proof}

We now turn to the proof of Proposition~\ref{prop:probability_container}. 

\begin{figure}[h!]
    \centering
    \includegraphics{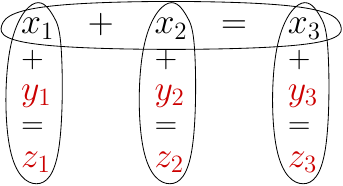}
    \caption{A wicket where all the $y_i$ and $z_i$ are red.}
    \label{fig:wicket_red}
\end{figure}

\begin{proof}[Proof of Proposition~\ref{prop:probability_container}]
We choose $\eps = \eps_{\ref{lem:lotsmc}}$ to be the constant given by Lemma \ref{lem:lotsmc}, which will be stated later, and fix an arbitrary container $C\in \mathcal{C}$.
We partition the elements of $[n]$ into four different sets: let $M_C$ be the elements \emph{missing} from the container $C$ (that means not being present in either red or blue), $R_C$ be the elements that are present \emph{only} in $C$ in the red copy of $[n]$, $B_C$ be the elements present only in the blue copy of $[n]$, and let $T_C$ be the \emph{two-coloured} elements (that is, those that are present in the container in both colours --- think of these as the elements where the container does not restrict the colouring).

Before diving into the details of the proof, we make some observations on what it means for a set $P$ to be compatible with the container $C$, thereby sketching our proof strategy. First, note that if $P$ contains an element that is missing from $C$, then clearly there is no colouring $\varphi$ of $A \cup P$, let alone a Schur colouring, such that $\varphi \subseteq C$, and so we do not have $P \triangleleft C$. In other words, the event that $P \triangleleft C$ implies that $P \cap M_C = \emptyset$.

Similarly, we can use the fixed red and blue elements of $C$ to derive restrictions on $P$ in the event that $P\triangleleft C$. Indeed, suppose that either $P\cap R_C$ or $P\cap B_C$ contains a Schur triple. Then, as the colour of these elements is predetermined by the container $C$, any colouring $\varphi$ of $A\cup P$ such that $\varphi\subseteq C$ will have a monochromatic Schur triple, and hence will not be a Schur colouring. Therefore, if $P\triangleleft C$, then both $P\cap R_C$ and $P\cap B_C$ must be sum-free. 

These simple implications will already allow us to handle some types of containers. Indeed, if a container $C$ is such that the set $M_C$ is linearly large, then it is highly unlikely that the random perturbation $P$ avoids it. Similarly, if $R_C$ or $B_C$ contains quadratically many Schur triples, then $P$ will almost surely contain one of them.

This leaves us with those containers for which $M_C$ is small and there are few Schur triples in $R_C$ and $B_C$, and this final case is more subtle. Using these conditions, together with the fact that $C$ spans few edges of the hypergraph $\mc H_A$, we will show that there are many wickets (recall Definition~\ref{def:wicket}) where the elements $y_i, z_i$ all have the same colour (say red); that is, they belong to $R_C$ (see Lemma~\ref{lem:lotsmc} for the precise statement). Then, with high probability, such a wicket appears in $P$, and this prohibits a Schur colouring of $A \cup P$. Indeed, if any of the elements $x_i$ are coloured red, they form a red Schur triple with $y_i$ and $z_i$. Otherwise, all the $x_i$ are blue, forming a blue Schur triple.

We now proceed to give the details of each case. 

\paragraph{Case I: $M_{C}$ contains at least $\eps n$ elements}
As discussed in the proof sketch above, the event that $P\triangleleft C$ is contained in the event that $P\cap M_C=\emptyset$. Hence, if there are at least $\eps n$ missing elements, we have that  \[\Pr[P\triangleleft C] \leq \Pr[P\cap M_C=\emptyset ]\leq (1-p)^{\eps n}\le e^{- \eps np }.\]
Thus, when $p = \omega \left( (sn)^{-1/3} \log n \right)$, we have the required bound \[\Pr[P\triangleleft C]\leq e^{-\omega\left(s^{-1/3}n^{2/3}\log n\right)}.\] The bound on $p$ holds due to the fact that $p=\omega\left((sn^{13})^{-1/27}\log n\right)$ and $s=\Omega(n^{1/2})$. 

\paragraph{Case II: either $R_{C}$ or $B_{C}$ contains at least $\eps n^2$ Schur triples}
Assume without loss of generality that $R_{C}$ has $\eps n^2 $ Schur triples and let $RT(C)$ be the set of all \emph{non-degenerate} Schur triples in $R_{C}$. As there are only $n/2$ degenerate Schur triples in $[n]$ and we are only interested in asymptotics, we may assume that $|RT(C)|\geq \tfrac{\eps}{2} n^2$. 
Now, as discussed above, if $P\triangleleft C$, then we must have that $RT(C)\cap P^3=\emptyset$. Hence, appealing to Lemma~\ref{lem:JansonSTs} with $\xi = \eps/2$, we conclude that the probability that $P\triangleleft C$ is at most $e^{-\Theta( n^2p^3)}$, and so we have the desired bound on $\Pr[P\triangleleft C]$ as long as $p=\omega((sn^4)^{-1/9}(\log n)^{1/3})$. The latter holds as $p=\omega\left((sn^{13})^{-1/27}\log n\right)$ and $s=\Omega(n^{1/2})$.

\paragraph{Case III: all remaining containers}
If a container $C$ is not covered by the previous two cases, then $|M_{C}|<\eps n$ and both $R_{C}$ and $B_{C}$ contain fewer than $\eps n^2$ Schur triples. The following lemma, which we shall prove later, shows that any such container must contain many wickets where, as illustrated in Figure~\ref{fig:wicket_red}, the elements $y_i$ and $z_i$ are all prescribed by the container to have the same colour.

\begin{lem} \label{lem:lotsmc}
There is some $\eps > 0$ such that applying Corollary~\ref{cor:containers} to $\mc H_A$ with $\eps$ yields the following. Let $C \in \mc C$ be a container for $\mc H_A$ for which $\card{M_C} < \eps n$ and $R_C$ and $B_C$ each contain fewer than $\eps n^2$ Schur triples. Then there is a constant $\xi > 0$ and a set $\chi \in \{ R_C, B_C \}$ such that there are $\xi n^5$ wickets where the elements $y_i, z_i$ belong to $\chi$.
\end{lem}

Assuming the lemma, the proof of Proposition~\ref{prop:probability_container} in this case is now straightforward. Indeed, we may without loss of generality take $\chi = R_C$. Suppose such a wicket was contained in $P$. Then, in any colouring contained in $C$, no element $x_i$ can be coloured red, as that would create a red Schur triple $(x_i, y_i, z_i)$. However, colouring each $x_i$ blue instead creates the blue Schur triple $(x_1, x_2, x_3)$. Hence, if $P \triangleleft C$, then $P$ cannot contain any of the $\xi n^5$ wickets given by Lemma~\ref{lem:lotsmc}. By Lemma~\ref{lem:JansonWickets}, this occurs with probability at most $e^{-\zeta n^5 p^9}$ for some constant $\zeta > 0$. By our choice of $p$, this is $e^{-\omega( s^{-1/3} n^{2/3} \log n) }$, as required.
\end{proof}

To prove Lemma~\ref{lem:lotsmc}, we will make use of the following claim, which asserts that we can find an interval on which one of the two monochromatic sets is quite dense.

\begin{clm}\label{clm:mcdense}
For $\eps \le 10^{-7}$, let $C \in \mc C$ be a container for $\mc H_A$ for which $\card{M_C} < \eps n$ and $R_C$ and $B_C$ each contain fewer than $\eps n^2$ Schur triples. Then there exists an $\eta \in [n]$ with $\eta \ge \tfrac12 n$ and a set $\chi \in \{ R_C, B_C \}$ such that $\card{ \chi \cap [\eta] } \ge \tfrac{9}{20} \eta$.
\end{clm}

We shall prove Claim~\ref{clm:mcdense} in due course, but let us first derive Lemma~\ref{lem:lotsmc} from it.

\begin{proof}[Proof of Lemma~\ref{lem:lotsmc}]
Choose $\eps = \tfrac{1}{4} \min\left\{\delta_{\ref{thm:removal}}\left(\tfrac{1}{100}\right), 10^{-8}\right\}$, where $\delta_{\ref{thm:removal}}\left(\tfrac{1}{100}\right)$ is the constant from Theorem~\ref{thm:removal} (Green's removal lemma) applied with $\eps_{\ref{thm:removal}} = \tfrac{1}{100}$, and let $\eta$ and $\chi$ be as given by Claim~\ref{clm:mcdense}, assuming without loss of generality that $\chi=R_C$. By assumption, $R_C \cap [\eta]$ has fewer than $\eps n^2 \le 4 \eps \eta^2$ Schur triples. By our choice of $\eps$, we can apply Theorem~\ref{thm:removal} to obtain a sum-free set $R'\subseteq R_C \cap [\eta]$ such that $\card{(R_C\cap[\eta])\setminus R'}\leq \tfrac{1}{100}\eta $. 
Hence, we have $|R'| \ge \card{R_C \cap [\eta]} -  \tfrac{1}{100}\eta > \tfrac{2}{5}\eta+1$, and Theorem~\ref{thm:sumstab} then yields that $R'$ is either contained in the odd numbers or does not have small elements (that is, $\min R' > \card{R'} > \tfrac{2}{5}\eta+1$).

We now construct the desired wickets. First, let $X$ be the set of all even integers not larger than $\tfrac{1}{10}\eta$. Given a pair of elements in $X$, their difference also lies in $X$. Thus, discounting the $\tfrac{1}{40}\eta$ Schur triples of the form $(x, x, 2x)$, we find at least $10^{-3} \eta^2$ triples $(x_1, x_2, x_3) \in X^3$ with $x_1 \neq x_2$ and $x_3 = x_1 + x_2$.

We next show that for any element $x\in X$, there are at least $\tfrac{1}{15}\eta$ Schur triples of the form $x+y=z$ with $y$ and $z$ distinct elements in $R'$. Then, to build one of the desired wickets, we can first choose a non-degenerate Schur triple $x_1 + x_2 = x_3$ in $X$, and subsequently choose for each $i$ a Schur triple $x_i + y_i = z_i$ with $y_i, z_i \in R'$. We shall need to ensure that these triples all use distinct elements (to obtain the nine-element wicket), but such considerations rule out at most a constant number of triples at each stage, and so we can build at least $10^{-8} \eta^5$ wickets in this fashion. Since $\eta \ge \tfrac12 n$, we can safely take $\xi = 10^{-10}$.

Given $x\in X$, suppose first that $R'$ is contained in the odd numbers. There are at least $\ceil{\tfrac12 (\eta - x)}$ odd numbers $y$ such that $z = x+y \le \eta$. Moreover, each odd number in $[\eta]$ appears in at most two of these triples. Since $\card{R'} > \tfrac25 \eta+1$, there are at most $\tfrac{1}{10} \eta$ odd numbers in $[\eta]$ missing from $R'$, and hence for at least $\tfrac12 (\eta - x) - \tfrac{1}{5} \eta$ of these pairs, we have $y,z \in R'$. As $x \in X$, we have $x \le \tfrac{1}{10} \eta$, and hence we have at least $\tfrac{1}{5} \eta$ Schur triples of the desired form in this case.

In the other case, $R'$ is such that $\min R' \ge |R'| > \tfrac25 \eta+1$, and we denote by $I$ the interval $\left[\ceil{\tfrac25 \eta}, \eta\right]$. There are $\floor{\tfrac35 \eta} - x$ integers $y \in I$ for which the sum $z = x + y$ lies in $I$ as well. Since $\card{R'} > \tfrac25 \eta+1$, there are at most $\floor{\tfrac15 \eta}$ integers in $I$ missing from $R'$, and each missing element appears in at most two of the triples. Hence, there are at least $\tfrac35 \eta-1 - x - \tfrac25 \eta$ pairs $y,z \in R'$ with $x + y = z$, and since $x \le \tfrac{1}{10}\eta$, this leaves us with at least $\tfrac{1}{15} \eta$ Schur triples, as required.
\end{proof}

The proof of Claim~\ref{clm:mcdense} is still outstanding, a situation we now rectify.

\begin{proof}[Proof of Claim~\ref{clm:mcdense}]
By Corollary~\ref{cor:containers}, we know that the container $C$ hosts at most $\eps sn^2$ edges of $\mc H_A$. Since each edge determines at most two targets in $A$, by averaging we can fix some element $\alpha\in A$ that is the target of at most $2\eps n^2$ edges. Our first aim is to find some $\eta\ge \tfrac{n}{2}$ and $Q \subseteq [\eta] \setminus \{\alpha\}$ such that the following hold:
\begin{align} \label{eq:largepartition}
    |Q| &\ge \tfrac{19}{20}\eta, \text{ and there is a partition } \Pi \text{ of } Q \text{ into pairs such that} \\
    \nonumber &\text{for every $\{\pi_1,\pi_2\}\in \Pi$, the set $\{\pi_1,\pi_2,\alpha\}$ hosts a Schur triple.}
\end{align} 
In other words, we aim to find some large $\eta$ such that almost all of the interval $[\eta]$ can be partitioned into pairs that form Schur triples with $\alpha$. We will then be able to show that there is some set $\chi\in \{R_C,B_C\}$ such that almost all the pairs of the partition $\Pi$ contain an element in $\chi$, which will complete the proof of the claim. 

Now, in order to establish \eqref{eq:largepartition}, we split into two cases, depending on how large $\alpha$ is. 
\paragraph{Case a: $\alpha\le \tfrac{1}{2}n$} 
Let $\ell = \left\lfloor\tfrac{n}{2\alpha}\right\rfloor $, set $\eta = 2 \alpha \ell$ and take $Q=[\eta]\setminus\{\alpha,2\alpha\}$.
We then choose 
\[ \Pi = \left\{\{2\alpha j +i, 2\alpha j + i + \alpha\}:0 \le j \le \ell-1, 1 \le i \le \alpha \right\}\setminus\{\{\alpha,2\alpha\}\}.\]
Clearly $Q$ and $\Pi$ have the required properties for \eqref{eq:largepartition} and $\eta \ge \max\{2\alpha,n-2\alpha\}\ge n/2$. 

\paragraph{Case b: $\tfrac{1}{2}n <\alpha$} Fix $\eta=\alpha$, $Q=[\eta]\setminus \{\floor{\tfrac{\alpha}{2}},\ceil{\tfrac{\alpha}{2}},\alpha\}$ and \[\Pi=\left\{\{i,\alpha-i\}:1 \le i \le \floor{\tfrac{\alpha}{2}} -1  \right\}.\]
Again, it is easy to check that $Q$ and $\Pi$ satisfy the conditions of \eqref{eq:largepartition}. 


\vspace{2mm}

Finally, given $\eta$, $Q$ and $\Pi$ as in \eqref{eq:largepartition}, we will show that there is a set $\chi\in\{R_C,B_C\}$ such that all but $2 \eps^{1/2} n$ pairs in $\Pi$ contain an element from $\chi$. Given the lower bound on the size of $Q$, the lower bound on $\eta$ and our upper bound on $\eps$, this will complete the proof of Claim~\ref{clm:mcdense}.

Firstly we define the following subsets of $\Pi$: 
\begin{align*}
    \Pi_0&=\{{\bf{\pi}}\in \Pi:{\bf{\pi}}\cap M_C\neq \emptyset\}, \\
     \Pi_R&=\{{\bf{\pi}}\in \Pi\setminus \Pi_0:{\bf{\pi}}\cap B_C= \emptyset\}, \\
        \Pi_B&=\{{\bf{\pi}}\in \Pi\setminus \Pi_0:{\bf{\pi}}\cap R_C= \emptyset\}, \\
           \Pi_{1}&=\{{\bf\pi}\in \Pi:|{\bf{\pi}}\cap B_C|=|{\bf{\pi}}\cap R_C|=1 \}.
\end{align*}
Note that $\Pi\subseteq \Pi_0 \cup \Pi_R \cup \Pi_B \cup \Pi_1$, but this is not quite a partition, as $\Pi_R$ and $\Pi_B$ intersect in pairs whose elements are both in $T_C$. By assumption, $|M_C|\le \eps n$, and so $|\Pi_0|\le \eps n$. Also, we have that either $\Pi_R$ or $\Pi_B$ contains at most $(2 \eps)^{1/2} n$ pairs. Indeed, observe that if $\pi_R = \{x,y\} \in \Pi_R$ and $\pi_B = \{z,w\} \in \Pi_B$, then the set $\{x_R, y_R, z_B, w_B\}$ forms an edge of $\mc H_A$ with target $\alpha$. By our choice of $\alpha$, there are at most $2 \eps n^2$ such edges, and so the smaller of $\Pi_R$ and $\Pi_B$ can contain at most $(2 \eps)^{1/2} n$ pairs.

Hence, without loss of generality, $|\Pi_B|\le (2\eps)^{1/2}n$. Since all pairs in $\Pi_1$ and $\Pi_R \setminus \Pi_B$ contain at least one element of $R_C$, and we can partition $\Pi$ as $\Pi_0 \cup \Pi_B \cup (\Pi_R \setminus \Pi_B) \cup \Pi_1$, it follows that
\begin{align*}\card{ ( \Pi_R \setminus \Pi_B ) \cup \Pi_1 } &= \card{\Pi_R \setminus \Pi_B } + \card{\Pi_1}\\
&= \card{\Pi} - \card{\Pi_0} - \card{\Pi_B} \\
&\ge \card{\Pi} - \eps n - (2 \eps)^{1/2} n \\
&\ge \card{\Pi} - 2 \eps^{1/2} n,
\end{align*}
completing the proof of the claim, and thereby the proposition. 
\end{proof}

\section{Concluding Remarks} \label{sec:conc}
\begin{figure}[h!]
    \centering
    \includegraphics[scale = 0.9]{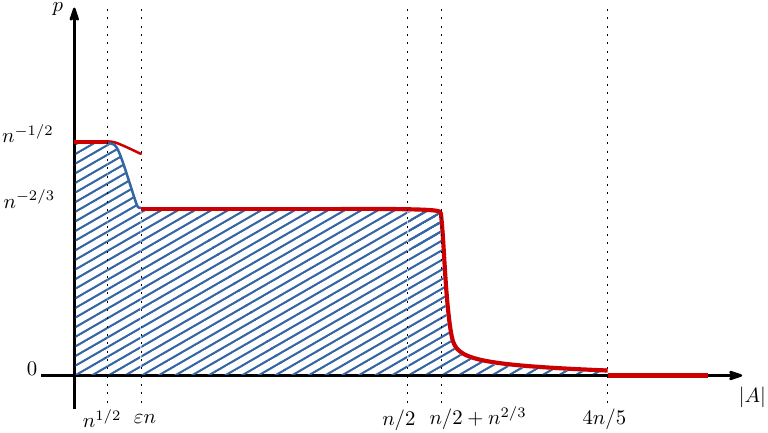}
    \caption{The complete picture so far}
    \label{fig:summary}
\end{figure}

In this paper, we explored Schur properties of randomly perturbed sets of integers. We addressed the case of sparse base sets $A$ and also dense base sets, describing the behaviour of the model as $|A|$ transitions from sublinear to linear and from $\tfrac{n}{2}$ to $\tfrac{4n}{5}$. A visualisation of our findings, as well as the previous works discussed in the introduction, can be seen in Figure~\ref{fig:summary}. We remark that our work completes the picture for dense base sets, bridging the gap between the extremal threshold of Hu (Theorem~\ref{thm:hu}) and the work of Aigner-Horev and Person (Theorem~\ref{thm:posdense}) giving a perturbed result that is tight for dense base sets $A$ with $|A|\le \tfrac{n}{2}$.

 On the other hand, for sparse base sets, we obtain upper and lower bounds that are a polynomial factor away from each other. 
We believe that the lower bound is likely to be the truth. Indeed, some evidence for this comes from our proof of the  $1$-statement of Theorem~\ref{thm:main_sparse}. There, we perform a union bound over choices of container that our random perturbation can be compatible with. Our analysis then splits into cases depending on properties of the container and in each case, we upper bound the probability of our random perturbation being compatible with such a container. Now in Case I, we look at containers that `miss' linearly many elements, forcing the random perturbation also to avoid these elements in order to be compatible. In this case our proof gives that a bound of $p=\omega((sn)^{-1/3}\log n)$ is already  sufficient to give our desired upper bound on the probability of compatibility. Therefore, if it were the case that all containers fall under Case I, then we would be able to provide a $1$-statement which gives the same bound as the $0$-statement up to a $\log$ factor. In fact, even more is true. If all containers were Case I, then we could remove the $\log$ factor appearing in the $1$-statement by  factoring in that for the random perturbation to be compatible with a container, it must also contain a small subset of the container (known as the `fingerprint' of the container, see the set $\mathcal{T}$ in Theorem~\ref{thm:container}). Similar ideas have been used in previous arguments using containers for sparse Ramsey theory (see for example~\cite{nenadov2016short}) and can be used to remove $\log$ factors in $1$-statement probabilities. 

Unfortunately, however, we see no reason why all containers should fall into Case I of our analysis. Indeed, one can cook up examples of candidate containers that satisfy property $2$ of Corollary~\ref{cor:containers} but do not miss any elements. For example, suppose $A$ is the set of the largest $s$ integers in $[n]$, and we create $C$ by taking all elements less than $\tfrac{n-s}{2}$ in both the red and blue copy of $[n]$ and all integers larger than $\tfrac{n-s}{2}$ in only the red copy $V_R$ of $[n]$. Such a set $C$ has no edges of $\mc{H}_A$ and does not fall into Case I (or Case II, for that matter). Despite the apparent necessity for more cases, a deeper analysis of the other cases could perhaps reduce the required probability for the $1$-statement, maybe all the way to match the $0$-statement. To summarise, we set the following problem.

\begin{prob}
Is it true that, for $s=s(n)\in \mathbb{N}$ and $p=p(n)$ such that $\Omega(\sqrt n)\le s\le n/2$ and $p=\omega((ns)^{-1/3})$, and for any $A\subseteq [n]$ with $|A|=s$,  we have that $A\cup [n]_p$ is Schur whp?
\end{prob}
Although we believe our argument falls short of the truth, we think that our approach for the $1$-statement is of value because it illustrates how the container method
can be combined with structural information about the underlying set of the hypergraph.
Note that improving the bound on the probability for Case III in the proof of Proposition~\ref{prop:probability_container} suffices to obtain a better $1$-statement.

\vspace{2mm}

Returning our focus to dense base sets, it would be interesting to extend the results to $r\ge3$ colours. As discussed in the introduction, whilst the random threshold is the same for all number of colours (see \cite{rodl1997rado}), the extremal threshold is already not known for $r\ge 3$. For the perturbed model, as noted by Aigner-Horev and Person~\cite{aigner2019monochromatic}, when $r\ge 3$, the $r$-Schur problem is only interesting for very dense base sets. Indeed, for $|A|\le n/2$, taking $A$ to be a sum-free set (and thus only using one colour for $A$) means that adding $o(n^{1/2})$ random elements does not help, as this random set can be $2$-coloured without a monochromatic Schur triple. On the other hand, with $\omega(n^{1/2})$ random elements, the random set is already $r$-Schur, without the need to consider the base set at all. 

Generalising this argument, we only see a separation in the behaviour of the randomly perturbed model and the purely random model for the Schur property with $r$ colours when the base set is large enough that it cannot be $(r-2)$-coloured without a monochromatic Schur triple.  Let $\mc{E}(r,n)$ be the extremal threshold for $r$ colours; that is, the minimum integer $m$ such that any subset $A\subset [n]$ of size at least $m$ is $r$-Schur. Then the following problem arises naturally.

\begin{prob} \label{prob:morecols}
For $r,n,s(n)\in \mathbb{N}$ with $r\ge 3$, determine  $p_r^*(n,s)$ such that the following statements hold. 
\begin{enumerate}
\item[(0)]  There exists a set $A \subseteq [n]$ with $\card{A}=\mc{E}(r-2,n) + s$ such that for  $p=o(p_r^*(n,s))$, whp $A\cup[n]_p$ is not Schur.
\item[(1)] For all $A \subseteq [n]$ with $\card{A}=\mc{E}(r-2,n) + s$ and $p=\omega(p_r^*(n,s))$, whp $A\cup[n]_p$ is Schur.

\end{enumerate}

\end{prob}

Note that as long as we can colour $A$ with $r-1$ colours without creating a monochromatic Schur triple, in order to obtain a set that is $r$-Schur, the perturbation probability needs to be at least $n^{-2/3}$, as otherwise the random set is sum-free whp. This shows that if $s$ is such that $\mc{E}(r-2,n) + s\le \mc{E}(r-1,n)$, we have $p_r^*(n,s)\geq n^{-2/3}$. It would be interesting to determine if the behaviour of $p_r^*$ is similar to what we observe here in the two colour case as the size of the base set moves beyond $\mc{E}(r-1,n)$. 


Whilst we believe it may be possible to make progress on Problem~\ref{prob:morecols} without knowing the values of the extremal thresholds, a better understanding of the extremal thresholds for $r\ge 3$ remains a central and very appealing problem in this area.


\begin{prob}
Determine $\mc{E}(r,n)$ for $r\ge 3$.
\end{prob}

\section*{Acknowledgements}
This project began at a workshop organised by Tibor Szab\'o from Freie Universit\"at Berlin, and the authors would like to thank him for his hospitality.

\bibliographystyle{elsarticle-num}
\bibliography{Biblio.bib}
\end{document}